\documentclass[11pt]{article}
\usepackage[english]{babel}

\usepackage{color}
\usepackage{array}
\usepackage{enumerate}
\usepackage{enumitem}
\usepackage{refcount}
\usepackage{algorithm}
\usepackage{algorithmic}
\usepackage{boxedminipage}
\usepackage{float}
\usepackage[utf8x]{inputenc}
\usepackage{soul}

\usepackage{color}
\usepackage{graphicx}
\usepackage{multirow}
\usepackage{tikz}
\usepackage{amsmath}
\usepackage{amssymb}
\usepackage{caption}
\usepackage{authblk}

\newtheorem{thm}{Theorem}[section]
\newtheorem{prop}[thm]{Proposition}
\newtheorem{lem}[thm]{Lemma}
\newtheorem{cor}[thm]{Corollary}
\newtheorem{defn}[thm]{Definition}

\newtheorem{exa}[thm]{Example}

\newtheorem{remark}[thm]{Remark}
\newtheorem{conj}[thm]{Conjecture}
\newtheorem{prob}[thm]{Problem}

\def\proof{{\bf Proof} }

\newcommand{\prfend}{\hbox to7pt{\hfil}
\par\vskip-\baselineskip\hbox to\hsize
{\hfil\vbox {\hrule width6pt height6pt}}\vskip\baselineskip}

\def\b{\par \noindent}

\def\dd{\medskip \par \noindent}
\long\def\eatit#1{}
\def\R{\mathbb{R}}

\def\N{\mathbb{N}}
\def\C{\mathbb{C}}
\def\A{\mathbb{A}}
\def\P{\mathbb{P}}

\def \"{``}

\font\tengothic=eufm10
\font\sevengothic=eufm7
\newfam\gothicfam
       \textfont\gothicfam=\tengothic
       \scriptfont\gothicfam=\sevengothic
\def\goth#1{{\fam\gothicfam #1}}

\usepackage{xypic}

\author[1]{Stefano Canino\thanks{stefano.canino@polito.it}}
\author[2]{Maria V. Catalisano \thanks{mariavirginia.catalisano@unige.it}}
\author[3]{Alessandro Gimigliano \thanks{alessandr.gimigliano@unibo.it}}
\author[3]{Monica Idà \thanks{monica.ida@unibo.it.}}
\affil[1]{Dipartimento di Matematica, Politenico di Torino, Corso Duca degli Abruzzi 24, 10129, Torino}
\affil[2]{Scuola Politecnica, Università di Genova, Via All'Opera Pia 15, 16145, Genova.}
\affil[3]{Dipartimento di Matematica, Università di Bologna, Piazza di Porta S.Donato 5, 40126, Bologna}

\begin{document}
\title{Superfat points and associated tensors.}
\maketitle

\begin{abstract}
We consider 0-dimensional schemes supported at a single point in $n$-space that are $m$-symmetric, i.e. that intersect any smooth curve passing through the point with length $m$, and the ones among them that are maximal with respect to inclusion (called $m$-superfat points). We study properties of such schemes, in particular for $n=2$. We give a first application of the simplest such schemes, namely 2-superfat points in the plane, by studying varieties defined by them on Veronese and Segre-Veronese varieties and the (symmetric or partially symmetric) tensors they parameterize.
 \end{abstract}
 
 \medskip
\noindent$Keywords:$
\smallskip 
\noindent 0-dimensional schemes; Secant varieties; W-states.

\medskip

\noindent {\small $MSC:$ \  Primary 14N07; 14M07. \ Secondary  14J81; 13H15.


\section*{Contents}

\quad\ \ 1. Introduction.

2. Symmetric and superfat points in $\P^n$.

3. Superfat and $m$-symmetric points in $\P^2$.

4. 2-squares on Veronese surfaces.

 \quad  4.1 The varieties $QQ(V_{d})$.  

5. 2-squares on Segre-Veronese surfaces.

 \quad  5.1 The first case: $q_2(V_{2,2})$. 
 
 \quad  5.2 The varieties $V_{d,d}$ and their $q_2(V_{d,d})$, $d\geq 3$.
 
 \quad  5.3 The varieties $qq_2(V_{d,d})$.
 
 \medskip


\section{Introduction}
The ideas for this work arose from problems encountered in a previous study, \cite{CGI}, related to the Jacobian \linebreak (or Tjurina) and Milnor schemes associated to a plane projective curve: in the case of curves with an ordinary singularity of multiplicity $m$, the Jacobian scheme at that point is $(m-1)$-symmetric, i.e. any line passing through its support intersects it with multiplicity $m-1$ (see Definition \ref{nmsymmetric}), while the Milnor scheme at that point is an $(m-1)$-superfat point (see Definition \ref{nsuperfat}; note that in \cite{CGI} such a scheme is called an \linebreak $(m-1)$-symmetric local complete intersection).  Inspired by these examples, we found the following questions quite natural to ask:
\begin{itemize}[leftmargin=*]
	\item What are the possible structures of a 0-dimensional scheme supported at one point $P\subset \P^n$ which are {\it $m$-symmetric}, i.e. that give the same length $m$ when intersected with any line passing through $P$?
	\item Given $m$, what is the structure of a maximal, with respect to inclusion, $m$-symmetric scheme, and what is the length of such a scheme?
\end{itemize}
In this paper we consider 0-dimensional schemes supported at one point $P$ in $\P^n_\C$ with respect to the properties related to the above questions. Given an $m$-symmetric scheme $X$ we show that if $X$ is maximal with respect to inclusion, then $X$ is a locally complete intersection of $n$ hypersurfaces with multiplicity $m$ at $P$, having tangent cones without common lines. In particular, it has length $m^n$. We call such a scheme an \"$m$-superfat point".
\dd 
We think that these schemes are interesting {\it per se}, e.g. for particular interpolation problems, but what we begin to study here is how to use these kind of schemes to parameterize structured tensors: when considering such a scheme $X$ in $\P^n$ or in a product of projective spaces, and then embedding the ambient space via Veronese or Segre-Veronese embeddings, we get that the linear span of the image of $X$ generates a linear space that lies between osculating spaces to the embedded variety, and such space parameterizes particular tensors. This aspect could be of interest also for possible applications to tensor decomposition.  In the secon part of the paper we begin to consider the simplest case: $m=n=2$, i.e. 2-symmetric schemes on Veronese or Segre-Veronese surfaces. 
\dd
The plan of the paper is the following: in Section 2 we give the main definitions and first properties of symmetric 0-dimensional schemes; in Section 3 we study more in detail the case of such schemes in the plane ponting out some non-intuitive properties of those schemes. Sections 4 and 5 are dedicated to the study of symmetric or partially symmetric $2^d$-tensors parameterized by points in the linear span of the image of 2-superfat points in $\P^2$ or in $\P^1\times \P^1$, via Veronese or Segre-Veronese embeddings, respectively. In particular, in section 4 we study the variety $Q(V_d)$, defined by the closure of the union of the linear spans of all the embeddings of 2-squares of $\P^2$ via its $d$-ple embedding; we show that $Q(V_d)=\tau_2(V_d)$, its second osculating variety, and we use this fact to show that $\tau_2(V_d)\subset \sigma_4(V_d)$, its 4-secant variety. Then we study the variety $QQ(V_d)$, parameterizing monomial forms in $\tau^2(V_d)$. In section 5 we consider particular 2-squares in $\P^1\times \P^1$, namely the ones defined by ideals $(\ell_1^2,\ell_2^2)$, where $\ell_1,\ell_2$ are lines in the two arrays of the product. We consider the  variety $q_2(V_{d,d})$ defined by the closure of the union of their images via the $(d,d)$-ple embedding of  $\P^1\times \P^1$, and we study its secant varieties and the subvariety $qq_2(V_{d,d})$ parameterizing monomial forms.

\section{Preliminaries on Symmetric and superfat points in $\P^n$}    
We will denote by $\P^n$ the $n$-dimensional projective space over the field $\C$ with coordinate ring $\C[x_0,\dots,x_n]$. Given a projective scheme $X\subset \P^n$, let  $I_X\subset\C[x_0,\dots,x_n]$ denote the defining ideal of $X$ and $\ell(X)$ its length. For any homogeneous ideal $I\subset\C[x_0,\dots,x_n]$, we denote by $I_d$ its homogeneous part of degree $d$.  In what follows, since the problems in consideration are actually of a local nature, when we work on a scheme $X$ supported on a point $P$ we usually suppose, without loss of generality, that $P=[1,0,\ldots ,0]$; in other words, we often work in the affine case, viewing $P=(0,\ldots ,0)$ in the affine chart $\{x_0\neq 0\}$, and using $(x_1,\ldots , x_n)$ as affine coordinates. 
\dd
We consider 0-dimensional schemes supported at one point $P$ which possess the following property.
\begin{defn}\label{nmsymmetric}\rm 
An  $m-$symmetric scheme is a 0-dimensional scheme $X$ supported  at one point $P \in \P^n$  and such that $\ell(X\cap L)=m$, for every line $L$ passing through $P$. 
\end{defn}
The easiest example of a scheme with this property is an $m$-fat point, i.e. a scheme $X$ such that $I_X = (I_P)^m$; such schemes are also denoted with $mP$, and we will use this notation. Let us recall that $\ell(mP)={m+n-1\choose n}$.  Actually it is immediate to check that the $m$-fat points are the $m$-symmetric schemes of smallest possible length, as shown by the following Lemma. 
\begin{lem}\label{lemmino} 
Let $X$ be an $m$-symmetric scheme  supported at $P\in \P^n$ and $I_X=(G_1,\dots,G_s)$ its defining ideal, where $\{G_1,\dots,G_s\}$ is a minimal set of generators. Then all the hypersurfaces $\{G_i=0\}$ have multiplicity at least $m$ at $P$, and at least $n$ of them have multiplicity exactly $m$. Moreover there is no line common to all the tangent cones of the hypersurfaces $\{G_i=0\}$ which have multiplicity exactly $m$ at $P$. 
\end{lem}
\proof
If there were a hypersurface $\{G_i=0\}$ having multiplicity at $P$ less then $m$ then a line $L$ passing through $P$ and not contained in the tangent cone of $\{G_i=0\}$,would locally intersect  $\{G_i=0\}$ with lenght  $m'< m$, hence we would have $ \ell(L\cap X)  \leq m' < m$, thus getting a contradiction.\\
Now let us suppose that there are $r$ hypersurfaces in $I_X$ with  multiplicity exactly $m$ at $P$, say \linebreak $\{G_1=0\},\dots,\{G_r=0\}$ and, by contradiction, that $r < n$. Let $\{F_1=0\}, \dots, \{F_r=0\}$  be their tangent cones at $P$ and consider the scheme $Y\subset \P^n$ defined by $I_Y=(F_1,\dots,F_{r})$. Since $Y$ is a cone and $\dim(Y)\geq 1$, then there is a line $L\subset Y$. By construction $L$ passes through $P$ and $\ell(L\cap X)>m$ against the hypothesis that $X$ is $m$-symmetric.\\
Finally if there were a line $L$ common to all the tangent cones $\{F_1=0\}, \dots, \{F_r=0\}$, then we would have $\ell(L\cap X)>m$, again contradicting our hypothesis of $m$-symmetry.  \prfend
\begin{lem}\label{lemmino2} 
Let $X$ be an $m$-symmetric scheme  supported at $P=(0,...,0)\in \A^n$. 
Then  there exists a regular sequence of $n$ forms of degree $m$ in the ideal of initial forms of $I_X\subset \C[x_1,\dots,x_n]$.
\end{lem} 
\proof  If we set $I_X = (G_1,...,G_s)\subset\C[x_1,\dots,x_n]$ then, by Lemma \ref{lemmino}, 
all $G_i$'s have multiplicity at least $m$ at $P$, and at least $n$ of them have multiplicity exactly $m$, so let, say, $G_1,\ldots,G_r$, $n\leq r\leq s$, be the ones that have tangent cone of degree $m$ at $P$, and let $F_1,\ldots,F_r \in \C[x_1,...,x_n]_m$ be the degree $m$ components of $ G_1, \ldots,G_r$, i.e. their tangent cones.  If the $F_i$’s are linearly dependent and, for example, $F_1= a_2F_2+ …. +a_rF_r$, then $G_1$ can be replaced by  $G_1-(a_2G_2+ …. +a_rG_r)$, which has tangent cone of degree $>m$; hence we can assume that the $F_i$’s are linearly independent, and, again by Lemma \ref{lemmino}, we have that $F_1, . . . , F_r$ do not have any common line.\\
We want to show that there are $n$ polynomials in $\langle F_1,\ldots F_r\rangle\subset\C[x_1,\dots,x_n]_m$ which have no common lines (actually we will find a regular sequence $H_1,\ldots , H_n$). This is obvious if $r=n$, so we can assume $r\geq n+1$. Consider $H_1=F_1$ and a generic linear combination of $F_1, . . . , F_{r}$:  $a_{11}F_1+a_{12}F_2+\ldots + a_{1r}F_{r} = H_2$. We want to check that $\dim (\{H_1=0\}\cap \{H_{2}=0\})=n-2$. 
Let $C_1,\ldots C_k$, $k\leq m$, be the irreducible components of $\{H_1=0\}$; let $P_i\in C_i \setminus P$, $\forall\; i=1,\ldots k$; in order to have that no $C_i$  is contained in $\{H_2=0\}$ is enough that $\forall\  i$, $H_2(P_i)\neq 0$, and this is true for the genericity of the linear combination $H_2$ (for each $i$, we can view  $x_{11}F_1(P_i)+x_{12}F_2(P_i)+\ldots + x_{1r}F_{r} (P_i)=0$ as a hyperplane in $\P^{r-1}$, with respect to homogeneous coordinates  $[x_{11},\ldots ,x_{1r}]$, so  it is enough to choose a point $[a_{11},a_{12}\ldots , a_{1r}]$ not lying on those hyperplanes).\\
Hence, $\{H_1=0\}\cap \{H_{2}=0\}$ has $\dim = n-2$. Now we repeat this procedure by definining  a generic linear combination $a_{21}F_1+a_{22}F_2+\ldots + a_{2r}F_r = H_3$ such that $\dim (\{H_1=0\}\cap \{H_{2}=0\}\cap\{H_3=0\}=n-3$ and so on, in order to get $H_1,\ldots , H_n$ that form a regular sequence and their intersection is only supported at $P$. \prfend
\begin{cor}\label{nmsymmetric2} Let $X$ be a 0-dimensional scheme supported at $P=(0,...,0)\in \A^n$, denote by \linebreak $I_X \subset \C[x_1,\ldots ,x_n]$ its ideal and set $\goth m=(x_1,\dots,x_n)$. The scheme $X$ is $m-$symmetric if and only if $I_X \subseteq \goth m ^m$ and the ideal of initial forms of $I_X$ contains a regular sequence of length $n$ of forms of degree $m$.
\end{cor}
\proof 
It is immediate to see that if $I_X \subseteq \goth m ^m$ and the ideal of initial forms of $I_X$ contains a regular sequence of length $n$ of forms of degree $m$, then $X$ is $m$-symmetric. Lemma \ref{lemmino2} gives the other implication.  \prfend
\begin{remark}\rm \label{fatinsym} As an immediate consequence of the previous lemma, we have that every $m$-symmetric scheme $X\subset \P^n$ supported at $P$ contains the $m$-fat point $mP$. 
Hence fat points are, with respect to inclusion, the smallest $m$-symmetric schemes; in particular, the length reaches its minimum, i.e., for every $m$-symmetric scheme $X$ we have $\ell (X) \geq {m+n-1 \choose n}$, with equality if and only if $X=mP$.
\end{remark}
Once seen which are the smallest $m$-symmetric schemes, we want to find out \"how big can an $m$-symmetric point be", i.e. we want to consider the following questions:
\begin{itemize}[leftmargin=*]
\item Among all the $m$-symmetric schemes supported on the same point $P$, which are the maximal ones with respect to schematic inclusion?
\item What is the maximum length of an $m$-symmetric scheme?
\end{itemize}
\begin{remark}\label{ci} \rm  The ideal $(x_1^{m},\ldots ,x_n^m)\subset\C[x_0,\dots,x_n]$ defines a projective scheme $X\subset\P^n$ of length\linebreak  $\ell(X)=m^n$, and it is easy to check that it satisfies $m$-symmetry, hence the answer to the second question above is at least $m^n$, i.e. the maximal length for an $m$-symmetric scheme in $\P^n$ is at least $m^n$. 
\end{remark} 
\begin{defn}\rm \label{nsuperfat}  A maximal, with respect to inclusion, $m$-symmetric scheme in $\P^n$ is called an {\it $m$-superfat point}, or just a superfat point if we do not need to specify $m$.
\end{defn}
\begin{thm}\label{nsuperfatmn}   An $m$-superfat point supported at $P\in \P^n$  is a locally complete intersection of $n$ hypersurfaces with multiplicity $m$ at $P$, having tangent cones without common lines. Thus, any $m-$superfat point in $\P^n$ has lenght $m^n$.
\end{thm}
\proof  
As usual, we can work in the affine case, with $P=(0,\ldots, 0)\in \A^n$. It is easy to see that a 0-dimensional scheme supported at $P$ which is locally complete intersection of $n$ hypersurfaces with multiplicity $m$ at $P$ and having tangent cones without common lines, is $m$-symmetric. Moreover, it has length $m^n$ by \cite{F2}, Corollary 12.4. Hence the statement is proved if we prove that any $m-$symmetric scheme $X$ supported at $P$  is contained in such a scheme. 
 \b In the notation of the proof of Lemma \ref{lemmino2}, 
let $K_1=G_1$ and $K_{i+1} =\sum_{j=1}^r a_{ij}G_j$, $i=1,\ldots ,n-1$, so $\{K_i=0\}$ has $\{H_i=0\}$ as tangent cone at $P$.  Since $\cap_{i=1,\dots,n} \{H_i=0\} =\{P\}$, the scheme $Z$ defined by the ideal $(K_1,\ldots , K_n)$ is 0-dimensional at $P$ and its component at $P$, $Z_P$, is locally complete intersection of $n$ hypersurfaces with multiplicity $m$, having tangent cones without common lines.  We have that \linebreak $(K_1,\ldots , K_n) \subset I_X$ and hence the schematic inclusion $X \subseteq Z_P$. \prfend
An immediate consequence is following corollary.
\begin{cor}\label{superfatmn2} Let $X$ be a 0-dimensional scheme supported at one point $P$ in $\P^n$. The following are equivalent:
\begin{enumerate}[leftmargin=*,label=\roman*.]
\item $X$ is an $m$-superfat point;
\item $X$ is a local complete intersection of $n$ hypersurfaces with multiplicity $m$ at $P$, having tangent cones without common lines;
\item $X$ is $m$-symmetric and $\ell (X)=m^n$.
\end{enumerate} 
\end{cor}
\proof One way follows from Theorem \ref{nsuperfatmn}. On the other hand, if (ii) holds, then $X$ has to be $m$-symmetric; by Lemma \ref{lemmino} its ideal is hence minimal among the ideals of all $m$-symmetric schemes supported on $P$, so $X$ is an $m$-superfat point.\prfend
\begin{remark}\rm Note that the Corollary above is not true if we work on the field $\R$, see below, Remark \ref{oddities}.
\end{remark}
We can check that being $m-$superfat can be described by the initial ideal of a scheme $X\subset \A^n$.
\begin{cor}\label{superfatalg} Let $X$ be a 0-dimensional scheme supported at $
P=(0,\ldots ,0)\in \A^n$, denote by \linebreak $I_X \subset \C[x_1,\ldots ,x_n]$ its ideal and set $\goth m=(x_1,\dots,x_n)$. The scheme $X$ is $m-$superfat if and only if $I_X \subseteq \goth m ^m$ and the ideal of initial forms of $I_X$ is generated by a regular sequence of length $n$ of forms of degree $m$.
\end{cor}
\proof If $X$ is a superfat point supported at $P$, Theorem \ref{nsuperfatmn} gives that $I_X \subseteq \goth m ^m$ and the ideal of initial forms of $I_X$ is generated by a regular sequence of length $n$ of forms of degree $m$.\\
Now assume that $X$ is a 0-dimensional scheme supported at $P$ with  $I_X \subseteq \goth m ^m$ and the ideal of initial forms of $I_X$ is generated by a regular sequence $g_1,\dots, g_n$ of forms of the same degree $m$. Note that this implies that $X$ is $m-$symmetric by Corollary \ref{nmsymmetric2}.
The $g_1,\dots,g_n$ are initial forms of $G_1,\dots, G_n \in I_X$, say $I_X = (G_1,\dots, G_n, \ldots , G_r)$.  
By Corollary \ref{superfatmn2} it is enough to show that $X$ is locally a complete intersection of $n$ hypersurfaces with multiplicity $m$ at $P$, having tangent cones without common lines. 
In $I_X$ we can suppose that any $G_i$ with $i>n$ has initial form of degree as large as we want, because, by adding combinations of the $G_1,\ldots G_n$, we can always increase the initial degree of such $G_i$, since the initial terms of $G_1,...,G_n$ do generate the ideal of initial forms of $I_X$, and this does not depend on the choice of its other generators. Let $k>0$ be such that $X \subseteq kP$ and assume that each $\{G_i=0\}$, $i>n$, has initial form of degree at least $k$; then, each $\{G_i=0\}$, $i>n$, contains $kP$ and hence $ X$. This implies that the scheme generated by $(G_1,...,G_n)$ is of type $X\cup X'$, where supp$(X'\cap P) = \emptyset$; so, locally, $X$ is generated by $(G_1,...,G_n)$ and we are done.
\prfend
Note that this shows how strong the hypothesis "The ideal of initial forms of $X$ is generated by a regular sequence" is. Following Remark 2.6, we define a more special class of superfat points.
\begin{defn}\rm \label{hypercube}  An $m$-superfat point whose ideal is of type $(l_1^m,l_2^m,\dots,l_n^m)$ for $l_i\in \C[x_0,...,x_n]_1$, with $l_1,\dots,l_n$ linearly indipendent,  is called an {\it $m$-hypercube}. If $n=2$, we use the notation ``$m-$squares", instead of $m-$hypercubes.
\end{defn}
In a sense, $m$-hypercubes are the simplest examples of $m$-superfat points.
\begin{remark} {\rm Two hypercubes of $\P^n$ given by the ideals $I=(l_1^m,...,l_n^m)$, $J=(h_1^m,...,h_n^m)$ with the same support (i.e. such that $(l_1,...l_n)= (h_1,...h_n)$), are different provided that $\{l_1,...l_n\}\neq \{h_1,...h_n\}$. In fact, $h_j^m \in I_m$ if and only if $\exists i$ such that $h_j=l_i$, since otherwise the forms $l_1^m,...l_n^m,h_j^m$, viewed as points of the Veronese variety $\nu_m(\P^n)$, are in general linear position. }
\end{remark}
\begin{remark} \rm  There are not only fat  and superfat points as hypercubes possessing $m$-symmetry; as an example, consider $X\subset \P^2$ defined by the ideal $(x^3,y^3,x^2y^2)$: this is $3$-symmetric and $\ell(X)=8$. 
\b A large class of examples is given by the Jacobian scheme of a plane curve at an ordinary singularity $P$ of multiplicity $m$, which is an $m-1$-symmetric scheme (see \cite{CGI}). We recall that (see \cite{BGM}, \cite{G}, \cite{LP}), if the curve has equation $f(x,y)=0$ in $\A^2$, the Jacobian scheme at $P$ is the component supported at $P$ of the scheme associated to the ideal $(f,f_x,f_y)$; its length $\tau$, called the Tjurina number, takes all the values in the interval 
$$\left\lfloor {3m^2-2m-4\over 4}\right\rfloor \leq \tau \leq (m-1)^2.$$
\b This gives the opportunity of pointing out a peculiar behaviours of the 0-dimensional schemes which sometimes baffle our intuition. 
Let $X=mP\subset \P^n$ and let $Y\subset \P^n$ be a different $m$-symmetric scheme in $\P^n$ supported at $P$ (which will have different length). We have that their linear sections $X\cap L = Y\cap L$ coincide for any line $L$ through $P$, but $X\neq Y$. If $J_L$ denotes the $m$-jet supported on $P$ and contained in $L$, we have  $Y\cap (\cup _{L\ni P} L)=Y\cap \P^n= Y$, and  $\cup_{L\ni P} (Y\cap L) = \cup_{L\ni P} J_L = mP=X$. Hence the schematic unions and intersections  $Y\cap (\cup _{L\ni P} L)$ and $\cup_{L\ni P} (Y\cap L)$ differs, while such unions and intersections are equal for sets.
\end{remark}
In defining $m-$symmetry we have used lines through the support point; the following result shows that this is equivalent to using smooth curves:
\begin{prop}  A 0-dimensional scheme $X$, supported at one point $P\in\P^n$, is $m$-symmetric if and only if $\ell(X\cap C)=m$ for every curve $C$ smooth at $P$.  
\end{prop}
\proof
Let $X\subset \P^n$ be an $m-$symmetric scheme with support at $P\in \P^n$, and $C\subset \P^n$ be a curve smooth at $P$.  We have $mP \subset X$, and $\ell(mP \cap C)= m$, so $\ell(X\cap C)\geq m$. Since $C$ is smooth, it is locally a complete intersection, i.e. there are polynomials $F_1,\ldots , F_{n-1}$ such that they are smooth at $P$, $(F_1,\ldots F_{n-1})$ defines $C$ at $P$, and the intersection of their tangent cones at $P$ is the tangent line $T_P(C)$. If we consider any $F\in I_X$, whose tangent cone has multiplicity $m$ and does not contain $T_P(C)$ (it must exist since $X$ is $m-$symmetric), then the length of the projective scheme defined by the ideal $(F,F_1,\dots,F_{n-1})$ is $m$ (by \cite{F2}, Corollary 12.4), and $(F,F_1,...,F_{n-1}) \subset I_X+I_C  = I_{X\cap C}$ hence $\ell(X\cap C) \leq m$. We are done, since the other implication is trivial.\prfend
Note that 0-dimensional schemes in $\P^n$ for $n\geq 3$ are not all smoothable (i.e. obtained by collapsing simple points); nevertheless, $m-$superfat points are all smoothable:
\begin{prop}\label{smoothable} Let $Z\subset \P^n$ be an $m-$superfat point. Then  $Z$ is smoothable, $\forall\  m,n\in \N$.
\end{prop}
\proof
The fact is actually well-known since every locally complete intersection 0-dimensional scheme is smoothable (e.g. see \cite{J}, Thm. 4.36) and, by Proposition \ref{nsuperfatmn}, $m$-superfat points are locally complete intersection. We just sketch the idea here: if we have, locally, $I_Z = (F_1,...,F_n)$, consider the schemes $Z_t$ defined, locally, by $I_{Z_t} = (F_1+tG_1,...,F_n+tG_n)$, where the $G_i$'s are generic forms of the same degree as  $F_i$. We have that, locally, $Z_t$ is given by $m^n$ simple points, and, as $t \rightarrow 0$, $Z_t \rightarrow Z$, so $Z$ is smoothable.\prfend
Now we want to check what happens when we consider the schematic union of all the $m$-hypercubes supported at a same point $P\in \P^n$.
\begin{thm}\label{union} Given a point $P\in \P^n$ and $m\geq 1$, we have that the schematic union of all $m$-hypercubes supported at $P$ is the fat point $(mn-n+1)P$.
\end{thm}
\proof
As usual, without any loss of generality, we can work in the affine case and consider the case \linebreak $P = (0,0,\ldots, 0)$.  Let $R:=\{(l_1,l_2,\ldots,l_n)\in\C[x_1,\ldots, x_n]_1^n\ |\ {\rm the}\ l_i's\ {\rm are\ linearly\ independent} \}$. Setting 
$$I := \bigcap _{(l_1,\dots l_n)\in R} (l_1^m,\ldots, l_n^m)$$
we have to prove that  $$I = (x_1,\ldots,x_n)^{nm-n+1}.$$
 First let us check that $ \displaystyle(x_1,\ldots,x_n)^{nm-n+1} \subset I$; actually, for any choice of $(l_1, \ldots,l_n)$ in $R$, every generator of  $(x_1,\ldots,x_n)^{nm-n+1}$ can be written as 
$$\sum_{i_1+\ldots+i_n=nm-n+1} a_{i_1,\ldots,i_n}l_1^{i_1}\ldots l_n^{i_n},$$
of course with different coefficients $a_{i_1,\ldots,i_n}\in \C$, if we change our choice of $(l_1,\ldots,l_n)$. Since \linebreak $n(m-1)= nm-n$, in every term of this polynomial at least one of the $l_i$'s appears with power $\geq m$; so we get $ \displaystyle(x_1,\ldots,x_n)^{nm-n+1} \subset I,$ as required.\\
Observe that $I$ is $GL_n$-invariant, and does not contain  $(x_1,\ldots,x_n)^{nm-n}$ (e.g. it does not contain the monomial $x_1^{m-1}\ldots x_n^{m-1}$); since the spaces $\C[x_1,\ldots,x_n]_k$ are the only $GL_n$-invariants inside $\C[x_1,\ldots,x_n]$, i.e. they are irreducible representations of $GL_n$, e.g. see \cite{FH}, we must have $I =  (x_1,\ldots,x_n)^{nm-n+1}.$ \prfend
Many examples let us think that the following conjecture is true.
\begin{conj} Given a point $P\in \P^n$ and $m\geq 1$, the schematic union of all $m$-superfat points supported at $P$ is the fat point $(mn-n+1)P$.
\end{conj}
Next we will consider the case $n=2$, where more detailed results are easier to get.
\section{Superfat and $m$-symmetric points in $\P^2$}When considering $m-$symmetric schemes and $m-$superfat points in $\P^2$, we use the coordinate ring $\C[x,y,z]$ rather then $\C[x_0,x_1,x_2]$, and often just $\C[x,y]$ when we work with ideals with support at one point, since we suppose it to be $P=[0,0,1]$ and we dehomogenize with respect to $z$. 
\dd
For $m=2$, we have that actually $2-$superfat points are $2-$squares:
\begin{prop}\label{2superfat2square} Every $2$-superfat scheme $Q\subset \P^2$ is a 2-square, i.e. $I_Q$ can be written (modulo projectivity) as $I_Q = (x^2,y^2)$.
\end{prop}
\proof
Let $P=[0,0,1]$ be the support of $Q$.  Since $\ell(Q)=4$,  there are at least 2 independent forms $F,G$ in $I_Q$ of degree 2. 
Since the ideal $(F,G)\subset I_Q$ defines a 0-dimensional scheme $X$ with $\ell(X) = 4=\ell(Q)$, we have  $(F,G)= I_Q$. The conics $\{F=0\}, \{G=0\}$ have a double point at $P$, otherwise their intersection would be a curvilinear scheme.
Let $F= L_1L_2$,  $G =L_3L_4$, where all $L_i\in \C[x,y]_1$. Now, in the pencil $\{aL_1L_2+bL_3L_4\}$ there will always be two conics of rank 1, since such pencil gives a line in $\P(\C[x,y]_2)\cong \P^2$ which will intersect in two points the conic representing the forms of rank 1 (i.e. the $2$-Veronese embedding of $\P^1$, parameterizing squares of linear forms). Note that our pencil cannot be represented by a tangent line to the conic since such lines represent pencils of conics with a common linear factor.  Hence every $2$-superfat point is always a 2-square and its ideal can be written, modulo projectivities, as $(x^2,y^2)$.
\prfend
\begin{remark} \rm This behavior (coincidence of 2-superfat points and 2-squares) is unique, i.e. there is nothing similar either in higher dimension or in higher degree in $\P^2$. E.g. in $\P^3$ the 2-superfat point of ideal $(xz^2, yz^2, z^3, xyz + z^2t, x^2, y^2)$ is not a 2-hypercube (it has generic Hilbert Function while a 2-hypercube does not). Actually the case $n=2=m$ is the only one when an $m-$hypercube has generic Hilbert function. The following example shows that for $m\geq 3$ in $\P^2$, $m-$squares and $m-$superfat points are distinct.
\end{remark}
\begin{exa}\rm Proposition \ref{2superfat2square} says that any 2-square is a global complete intersection in $\P^2$ of 2 conics. For $m=3$, the situation differs from the $m=2$ case, i.e. a $3$-superfat scheme $X$ (of length $9$) is not always a global complete intersection of 2 cubics.  Of course any ideal of type $I_X = (L_1L_2L_3,M_1M_2M_3)$, where $L_i,M_j\in \C[x,y]_1$ and $L_i\neq \alpha M_j$, for all $i,j\in \{1,2,3\}$, gives a c.i. $3$-symmetric scheme $X$ of length 9, but they are not all.  Consider the ideal $I_X=((x - y)^3, y^3z + x^2y^2,xy^3, y^4)$; it can be seen  that it defines a scheme of length 9 (e.g. one can check this by using CoCoA) which  is $3$-symmetric, i.e. a 3-superfat point. Nevertheless, $X$  is not a complete intersection, since its ideal generation is the generic one for a scheme of length 9: one cubic and three quartics. Anyway, if we consider the intersection of its two first generators at $P=[0,0,1]$, we get a scheme of length 9, which has to be $X$; in other words, the scheme $X$ is  the local complete intersection of two curves with a triple point at $P$ and with no common tangent, in accord with Theorem \ref{nsuperfatmn}.
\end{exa}
Actually, we can have $m$-superfat points in $\P^2$ with generic Hilbert function only for $m\leq 3$, as shown by the following proposition.
\begin{prop}\label{maxHilb} For every $m\geq 4$, there does not exist an $m-$superfat point in $\P^2$ having maximal Hilbert function.
\end{prop}
\proof
Let $X\subset \P^2$ be an $m-$superfat point with support at $P$. For $m\geq 4$, we have $m^2>{m+2\choose 2}$, hence if $X$ has generic Hilbert function its ideal is generated in degrees $\geq m+1$, so $I_X\subset I_P^{m+1}$, i.e. $(m+1)P \subset X$, which is absurd since $X$ is $m-$symmetric.
\prfend
\begin{remark}\label{oddities}\rm Recall that by Theorem \ref{union}, we have that the union of all $2-$squares with the same support $P$ is the scheme $3P$. Let us note that something quite different can happen if we do not consider all the couple of lines as we did in Theorem \ref{union}. For example, consider $P=[0,0,1]\in \P^2$ and the union of the 2-squares supported at $P$ that are defined via two lines which are \"perpendicular" with respect to the apolar action of $\C[x,y]$ on itself (i.e. when we view $x,y$ as the derivations $\frac{\partial}{\partial x}$,$\frac{\partial}{\partial y})$; in this case we do not get the entire fat point $3P$. 
\b To say that $\ell_1,\ell_2$ are linear perpendicular forms means that if $\ell_1 = ax-by$, then $\ell_2 = bx+ay$; if we moreover ask that $(\ell_1^2,\ell_2^2)$ is a 2-square, the lines defined by $\ell_1,\ell_2$ are different. Hence, if we set \linebreak
$\mathcal{S} = \{(\ell_1,\ell_2)\in (I_P)_1\times (I_P)_1 |\ \ell_1\perp\ell_2 , (\ell_1^2,\ell_2^2)\ is\ a\ 2-square \}$ we get  $$\bigcap _{(\ell_1, \ell_2) \in \mathcal{S} }( \ell_1^2,\ell_2^2) = (x^2+y^2, x^3, x^2y).$$  
In fact, it is quite immediate that each ideal $(\ell_1^2,\ell_2^2)$ contains the ideal $(\ell_1,\ell_2)^3 =(x,y)^3$; if moreover  $\ell_1\perp\ell_2$, we can write $\ell_1 = ax-by$, $\ell_2 = bx+ay$, hence $(\ell_1^2,\ell_2^2)$ contains both $a^2x^2+b^2y^2-2abxy$, $b^2x^2+a^2y^2+2abxy$, and so  $(a^2+b^2)(x^2+y^2)\in (\ell_1^2,\ell_2^2)$, i.e. $(x^2+y^2)$ is contained in any ideal $(\ell_1^2,\ell_2^2)$ with $(\ell_1,\ell_2)\in\mathcal{S}$, and the thesis follows.  
\dd Note that it is actually enough to intersect two of those ideals to obtain the total intersection ideal. Note also that $(a^2+b^2)\neq 0$, since the pairs $(a,b) = \alpha(1,i)$ for which it is zero correspond  to the only two particular lines through $P$, namely $\{x\pm iy =0\}$ which we have to exclude among the pairs of lines in $\mathcal{S}$ , because they are isotropic, i.e. \"perpendicular to themselves", hence the ideal $((x\pm iy)^2,(ix\mp y)^2) = ((x\pm iy)^2)$ is not the ideal of a $2-$square but of a double line. 

\dd It is also interesting to observe that the scheme  $Z = \bigcup _{(\ell_1, \ell_2) \in \mathcal{S} } Q_{\ell_1\ell_2}$, where $I_{Q_{\ell_1\ell_2}}= ( \ell_1^2,\ell_2^2)$, is not \linebreak 2-symmetric, even if it is an (infinite) union of 2-symmetric schemes: $Z\cap L$ has length 2 for all lines $L$, except for the two lines $x\pm iy=0$ which meet it with length 3.\\
On the other hand, the fat point $3P$ is 3 -and not 2- symmetric, although it is union of 2-symmetric schemes by Theorem \ref{union}.
\dd
Finally, observe that the scheme $Z$ considered above {\em is} 2-symmetric if we consider it over the reals. Hence Theorem \ref{nsuperfatmn} does not hold over $\mathbb{R}$, in fact $(x^2+y^2,x^2y,x^3)$ defines a 2-symmetric scheme in $\P^2_\R$ of length 5. 
\end{remark}
\begin{remark}\label{post}\rm In a forthcoming paper the authors, together with Alessandro Oneto, prove that a generic union of $s$  2-squares (where generic means that the position of the $s$ support points is generic and so are the 2s lines defining the ideals of the 2-squares) has good postulation, i.e. has maximal Hilbert Function.
\end{remark}
\section{$2$-squares on Veronese surfaces}
Now we want to begin to see how the square points can give, with their immersions on Veronese surfaces, parameterizations of structured tensors. We will considering here only 2-squares, which are the ones with the  \"best behaviour".
\dd
Let us recall the following definitions which we are using here and in the following section.   
\begin{defn}\label{tauandsigma}\rm Osculating and secant varieties. Let $X\subset \P^n$ be a smooth variety.
\begin{itemize}[leftmargin=*]
\item  The $k^{th}$-$osculating$ $variety$ $\tau_k(X)$ is $\tau_k(X) := \overline{\bigcup_{P\in X} \tau_{k,P}(X)}$, where $\tau_{k,P}(X) = \langle X \cap (k+1)P\rangle$, is the span of the $k^{th}$ infinitesimal neighborhood of $P$ on $X$ (so $\tau_0(X) = X$ and $\tau_1(X)$ is its tangential variety); 
\item Let us denote with $\mathbf{P}_k = \{P_1,\ldots,P_k\}\subset X$ any set of distinct points, and $\sigma^0_k(X) := \bigcup_{\mathbf{P}_k\subset X} <\mathbf{P}_k>$. Then $\sigma_k(X) := \overline{\sigma^0_k(X)}$ is the $k^{th}$-$secant$ $variety$ of $X$ (the variety of $(k-1)$-subspaces $k$-secant to $X$); 
\item for a point $T \in \P^n$, we say that the $X$-$rank$ of $T$ is $r$, and we write $rk_X(T)=r$, if $r = \min\{ k| T \in \sigma^0_k(X)\}$, while the $X$-$border$ $rank$  of $T$ is $r$, and we write $bdrk_X(T)=r$, if $r = \min\{ k| T \in \sigma_k(X)\}$.
\end{itemize}
\end{defn}
We will consider the varieties $V_{n,d}$, the $d^{th}$ Veronese embedding $\nu_{n,d}$ of $\P^n$ into $\P^{N_{n,d}}$, $N_{n,d} = {n+d\choose n}-1$. We view $\P^n$ as $\P(\C[x_0,\ldots x_n]_1)$ with homogeneous coordinates $[y_0,\ldots ,y_n]$, so that any $P= [y_0,\ldots ,y_n]\in \P^n$ corresponds to a linear form $L_P = y_0x_0+\ldots + y_nx_n$ and $\nu_{n,d}(P)$ to the form $$L_P^d = \sum_{\mathbf{i}\in \N^{n+1}; |\mathbf{i}|=d} {d\choose \mathbf{i}}y_0^{i_0}\ldots y_n^{i_n}x_0^{i_0}\ldots x_n^{i_n},$$
where ${d\choose \mathbf{i}} := \frac{d!}{i_0!i_1!\ldots i_n!} =  {d \choose i_0}{d-i_0\choose i_1 }\cdots{d-i_0-\ldots i_{n-1}\choose i_n}$. \\
What we mean is that $\nu_{n,d}(P) = (y_0^d,\ldots ,{d\choose \mathbf{i}}y_0^{i_0}\ldots y_n^{i_n},\ldots ,y_n^d)\in \P^{N_{n,d}}$, endowed with homogeneous coordinates $[z_{i_0,\ldots,i_n}]$, $\mathbf{\alpha} =(i_0,\ldots i_n)\in \N^{n+1}$ with $|\mathbf{\alpha}| =d$, ordered by the usual lexicografic order, so that each point $[z_{d,0,\ldots,0},\ldots,z_{0,\ldots ,0,d}]$ parameterizes the form $\sum_{\mathbf{i}\in \N^{n+1}; |\mathbf{i}|=d} {d\choose \mathbf{i}}y_0^{i_0}\ldots y_n^{i_n}z_{i_0,\ldots, i_n}$.  In this way, $V_{n,d}$ parameterizes the forms of degree $d$ which can be written as a $d^{th}$-power of a linear form. In the language of tensors, $\P^{N_{n,d}}$ parameterizes symmetric $(n+1)^d$-tensors (i.e. ${\underbrace{(n+1)\times \ldots\times (n+1)}_{d\ times}}$ symmetric tensors, and $V_{n,d}$ the decomposable ones. 
\begin{defn}\label{rkandbdrk} \rm  
 For any form $F \in \C[x_0, \ldots , x_n]_d$, consider the corresponding point $T_F\in \P^{N_{n,d}}$, i.e. the point parameterizing $F$; we say that the (Waring) rank of $F$ is $r$ (or $rk(T_F) = r$) if $r = rk_{V_{n,d}}(T)$. We have the analogous definition for the border rank $bdrk(T_F)$.
\end{defn}
Note that when we view the form $F$ above as associated to a symmetric tensor $T$, this correspond also to the notion of \"symmetric rank" ($srk(T)$) of the tensor.
\dd
Here we consider the Veronese embeddings of the plane $\P^2$, i.e. the surfaces $V_{2,d}$, so we will write simply $V_{d}:=V_{2,d}$, $\nu_d := \nu_{n,d}$ and $N_d := N_{2,d}$. 
If we have a $2$-square $Q_P\subset \P^2$, then $<\nu_d(Q_P)>\cong \P^3$. Actually, we can identify $<\nu_d(Q_P)>$ with $(I^\perp_{(Q_P)})_d$, where the $\perp$ is considered with respect to the apolarity action (see e.g. \cite{BCCGO}). More specifically, let us consider $\mathbf{\alpha} =(i,j,k)\in \N^{3}$ and $\nu_d([x_0,x_1,x_2]) = [x_0^d,\ldots ,{d\choose \mathbf{i}}\mathbf{x}^\mathbf{\alpha},\ldots ,x_2^d]$ (where $\mathbf{x}^\mathbf{\alpha} = x_0^{i}x_1^{j}x_2^{k}$), and homogeneous coordinates $[z_{d,0,0},z_{d-1,1,0},z_{d-2,2,0},...,z_{0,0,d}]$ in $\P^{N_d}$, so that $V_{d}$ has parametric equations $z_{i,j,k} = \mathbf{x}^\mathbf{\alpha} = x_0^ix_1^jx_2^k$, for all $(i,j,k)\in \N^3$ with $i+j+k=d$. If we choose coordinates so to have $I_{Q_P}=(x_0^2,x_1^2)$, we can check that $<\nu_d(Q_P)> = (I_{Q_P}^\perp)_d \cong <x_2^d, x_2^{d-1}x_0,x_2^{d-1}x_1, x_2^{d-2}x_0x_1>$. Actually, here $<\nu_d(Q_P)>$ is defined by equations $z_{i,j,k}= 0 ,\  \forall\; k\leq d-3,$ and $z_{i,j,d-2}=0,\  \forall\; (i,j)\neq (1,1)$.
\dd
Now we want to consider the variety spanned by all the possible schemes $\nu_d(Q_P)$, on the surface $V_{d}$. We have the following:
\begin{prop}\label{Formequadrellate} Let  $Q(V_d) := \overline{\bigcup_{Q_P\subset \P^2} <\nu_d(Q_P)>}$; then we have  $Q(V_d) = \tau_2(V_d)$, the 2-osculating variety to $V_d$.  Moreover, if a point $T$ of $P^{N_d}$ lies on $\tau_2(V_d)$, then it parameterizes a form  $F_T\in \C[x_0,x_1,x_2]_d$ which, modulo a change of variables in $\P^2$, can be written either as $x_2^{d-2}(a_2x_2^2+a_1x_2x_1+a_0x_0x_2+a_3x_0x_1)$ or as $x_2^{d-2}(a_1x_1x_2 + a_0x_0^2)$. 
\end{prop}
\proof  Actually, the union of all $2$-squares supported at the same point in $\P^2$ is the fat point $3P$, see Proposition \ref{union}, and $<\nu_d(3P)>$ spans the $\P^5$ osculating $V_d$ at $\nu_d(P)$, hence  $Q(V_d)\subset \tau_2(V_d)$. 
Moreover, every  $T\in \tau_2(V_d)$ parameterizes a form which can be written as $F_T=x_2^{d-2}G$, where $G$ is a conic. \linebreak If $G = a_2x_2^2+a_1x_2x_1+a_2x_2x_0+ H(x_1,x_0)$, with $H\in \C[x_0,x_1]_2$, which is not a square, then \linebreak $H=(\alpha x_0+\beta x_1)(\gamma x_0+ \delta x_1) = \ell_0\ell_1$, with $\ell_0\neq \ell_1$, and we can write $x_2(a_0x_0+a_1x_1) = x_2(b_0\ell_0+b_1\ell_1)$, hence $F_T$ can be written as $x_2^{d-2}(a_2x_2^2+b_1x_2\ell_1+b_0x_2\ell_0+\ell_1\ell_0)$, as required (here $T\in <\nu_d(Q_P)>$, with $I_{Q_P} = (\ell_0^2,\ell_1^2)$ ).\\ 
The only exception is when $H = \ell_0^2$ for some $\ell_0\in \C[x_0,x_1]_1$. In this case we get 
$$F_T= x_2^{d-2}(x_2(a_2x_2+a_1x_1+a_0x_0) + \ell_0^2)$$
that, modulo a linear change of coordinate, we can write $x_2^{d-2}(b_2x_2x_1 +b_0x_0^2)$, as required in the statement. 
In this case  $T\notin \langle\nu_d(Q_P)\rangle$, for any $Q_P\subset \P^2$, nevertheless $T\in \overline{\bigcup_{Q_P\subset \P^2} <\nu_d(Q_P)>}$, in fact consider the linear form $\ell_\epsilon = x_0+\epsilon x_1$, and $Q_P = (\ell_\epsilon^2,x_0^2)$, then we have $\langle \nu_d(Q_P)\rangle = \langle x_2^{d}, x_2^{d-1}x_0, x_2^{d-1}\ell_\epsilon, x_2^{d-2}x_0\ell_\epsilon\rangle$, so $\lim_{\epsilon\rightarrow 0}\langle \nu_d(Q_P)\rangle = \langle x_2^{d}, 2x_2^{d-1}x_0, x_2^{d-2}x_0^2\rangle$ and any $T$ in there is such that $F_T$ can be written (modulo projectivities on the plane) as $x_2^{d-2}(a_1x_2x_1 + a_0x_0^2)$. Hence, also those $T$'s are in  $Q(V_d) := \overline{\bigcup_{Q_P\subset \P^2} <\nu_d(Q_P)>}$. \prfend
There is a noteworthy consequence of Prop. \ref{Formequadrellate} (for $d=4$ this is Lemma 41 in \cite{BGI}): 
\begin{cor}\label{Osculosecant} The second osculating variety $\tau_2(V_d)$ of a Veronese surface $V_d\subset \P^{N_d}$ is contained in the secant variety $\sigma_{4}(V_d)$. Thus, every $T\in \tau_2(V_d)$ has border rank $\leq 4$.
\end{cor}
\proof This is a direct consequence of the previous proposition, since  $\tau_2(V_d)=Q(V_d)$ and if $T\in <\nu_d(Q_P)>$, for some $Q_P$, then it is trivially contained in  $\sigma_{4}(V_d)$, since $\nu_d(Q_P)$ is a 0-dim scheme of length 4 and it is smoothable (all 0-dim schemes in $\P^2$ are), hence $<\nu_d(Q_P)>$ is the limit of a family of $\P^3$'s which are 4-secant to $V_d$. 
On the other hand, if $T$ is not on any $<\nu_d(Q_P)>$, by Prop. \ref{Formequadrellate}, let $F_T=x_2^{d-2}(a_1x_2x_1 + a_0x_0^2)$; if $a_1=0$, $F_T = x_2^{d-2}x_0^2$, then $T\in \tau_2(C_d)\subset \sigma_3(V_d)\subset \sigma_4(V_d)$ for some rational normal curve $C_d\subset V_d$. If $a_0= 0$, then $T\in \tau_1(V_d) \subset \sigma_2(V_2)$. When both $a_1,a_0$ are not zero, consider the scheme $Z\subset \P^2$, with \linebreak $I_Z = (x_1^2, x_0x_1, x_0^3)$. We have $\ell(Z) = 4$, and 
$$\langle \nu_d(Z)\rangle = (x_1^2, x_0x_1, x_0^3)^\perp _d = \langle x_2^d, x_2^{d-1}x_0, x_2^{d-1}x_1, x_2^{d-2}x_0^2\rangle \subset \sigma_4(V_d).$$
 So $F_T\in  \langle x_2^d, x_2^{d-1}x_0, x_2^{d-1}x_1, x_2^{d-2}x_0^2\rangle$, and also $T\in \sigma_4(V_d)$.  Note that 
 $$\langle x_2^d, x_2^{d-1}x_0, x_2^{d-1}x_1, x_2^{d-2}x_0^2\rangle = \langle x_2^d, x_2^{d-1}x_0, x_2^{d-1}x_1\rangle + \langle  x_2^d, x_2^{d-1}x_0, x_2^{d-2}x_0^2\rangle = \tau_{1,P}(V_d)+\tau_{2,P}(C_d),$$
  i.e. it is a $\P^3$, sum of two planes which have the line $\tau_{1,P}(C_d) =  \langle x_2^d, x_2^{d-1}x_0\rangle$ in common.
\prfend 
\begin{remark}\rm
The above results cannot be generalized to $3-$squares (or higher $m-$squares); e.g. if we were to consider the variety $Q_3(V_d) := \overline{\bigcup_{Q_{3,P}\subset \P^2} <\nu_d(Q_{3,P})>}$, $d\geq 7$, where the $Q_{3,P}$'s are the $3-$squares, we do not get $Q_3(V_d) = \tau_4(V_d)$, even if for each $P\in \P^2$, $\bigcup Q_{3,P} = 5P$; since $\ell(Q_{3,P}) = 9$, we have $Q_3(V_d) \subset \sigma_9(V,d)$, hence $\dim Q_3(V_d) \leq \dim \sigma_9(V_d) = 26$, while $\dim \tau_4(V_d)= 28$.
\end{remark}
Consider multi-indices  $(i_0,\ldots, i_n)\in \N^{n+1}$, such that $i_0+\ldots + i_n =d$, and recall (e.g. see \cite{Pu}) that for a positive integer $i<d$, the catalecticant matrix $Cat(i,d-i; n+1)$ is the ${n+i\choose n}\times {n+d-i\choose n}$-matrix with row and column indices respectively given by the multi-index sets (lexicografically ordered):  $\{\mathbf{\alpha} |\  |\mathbf{\alpha}| = i\}$ and $\{\mathbf{\beta} |\  |\mathbf{\beta}| = d-i\}$, where its $(\mathbf{\alpha}, \mathbf{\beta})$-entry is equal to $a_{\mathbf{\alpha}+\mathbf{\beta}}$. Note that we have $Cat(i,d-i; n+1) = {}^tCat(d-i,i; n+1)$.\\
When we have a form $F\in \C[x_0,\ldots x_n]_d$, its $(i,d-i)$-catalecticant matrix is the matrix\linebreak $M_i=Cat_F(i,d-i; n+1)$ such that
$$F =  (x_0^i,x_0^{i-1}x_1,\ldots x_n^i)\cdot M_i\cdot {}^t\!(x_0^{d-i},x_0^{d-i-1}x_1,\ldots , x_n^{d-i}).$$
If we view $F$ as a $(n+1)^d$-symmetric tensor, associated to a $(n+1)^d$-symmetric array, its $Cat_F(i,d-i; n+1)$ catalecticant matrix is obtained by flattening the array,
and eliminating equal rows or columns.\\
For every  $F\in \C[x_0,x_1,x_2]_d$, $d\geq 3$, its $(2,d-2)$-catalecticant matrix $M_2$ is of the form:
$$M_2 = \left( \begin{array} {cccccccc}  
a_{d,0,0}&a_{d-1,1,0}&\ldots&a_{2,0,d-2}\cr
 a_{d-1,1,0}&a_{d-1,0,1}&\ldots&a_{1,1,d-2}\cr 
 a_{d-1,0,1}&a_{d-2,2,0}&\ldots&a_{0,2,d-2}\cr
   a_{d-2,2,0}&a_{d-2,1,1}&\ldots&a_{1,0,d-1}\cr 
 a_{d-2,1,1}&a_{d-2,0,2}&\ldots&a_{0,1,d-1}\cr
  a_{d-2,0,2}&a_{d-3,3,0}&\ldots&a_{0,0,d}\cr
 \end{array}\right).$$
Now let us consider the $M_2$ catalecticant matrix associated to a generic form $F$ such that $T_F \in \tau_2(V_d)$; since by Prop. \ref{Formequadrellate} we can write (modulo projectivities in $\P^2$), $F =x_2^{d-2}(a_2x_2^2+a_1x_2x_1+a_0x_0x_2+a_3x_0x_1)$, $M_2$ is a $6\times {d\choose 2}$-matrix which  can be written in such a way that only the last five columns have some non-zero entries (for the case $d=4$, see also \cite{BGI}, Thm. 4.4 (2)):

$$M_2 = \left( \begin{array} {cccccccc}  0&\ldots&0&0&0&0&0&0\cr
 0&\ldots&0&0&0&0&0&a_{1,1,d-2}\cr 
 0&\ldots&0&0&0&0&a_{1,1,d-2}&0\cr
   0&\ldots&0&0&0&0&0&a_{1,0,d-1}\cr 
 0&\ldots&0&0&a_{1,1,d-2}&0&0&a_{0,1,d-1}\cr
  0&\ldots&0&a_{1,1,d-2}&0&a_{1,0,d-1}&a_{0,1,d-1}&a_{0,0,d}\cr
 \end{array}\right).$$
 By the way, there is a mistake in \cite{BGI}, where it is stated that such tensors can be written as $x_0^2x_1x_2$ (via a Gauss elimination on $M_2$); this is false since that Gauss elimination does not correspond to a projectivity in $\P^2$. We will analyze which tensors are of that monomial type in Prop.\ref{MonomialForms}.

\eatit{\medskip
The secant varieties to $\tau_2(V_d)$ have been studied, e.g. see \cite{BCGI}, \cite{BCGI2}, \cite{BF}), and the results there yield the following:

\bigskip
\begin{cor} Every form in $\C[x_0,x_1,x_2]_d$ can be written as a sum of $s= \lceil \frac{d^2+3d+2}{16}\rceil$, tensors of the form described in Prop. \ref{Formequadrellate} with the only exception of the case $d=4$, for which such $s$ is $3$ and not $2$.
\end{cor}
\proof  We have $dim (\tau_2(V_d)) = 7$, and from \cite{BF}, it is known that  $\sigma_s(\tau_2(V_d))$ has always the expected dimension, except for the case of $\sigma_2(\tau_2(V_4))\subset \P^{14}$, which should fill up its ambient space and it is a hypersurface, instead. This implies that, apart from the case $d=4$, for which we need $\sigma_3(\tau_2(V_4))$ to fill up $\P^{14}$, in all other cases $\dim( \sigma_s(\tau_2(V_d)))= 7s + s-1 = 8s-1$, hence we get that the first $s$ for which $\sigma_s(\tau_2(V_d))$ is the whole ambient space is for $8s-1\geq {d+2 \choose 2}-1$, i.e.  $s= \lceil \frac{d^2+3d+2}{16}\rceil$. \prfend }
\subsection{The varieties $QQ(V_{d})$}
The variety $\tau_2(V_d)$ contains a 1-codimensional subvariety parameterizing more particular forms, namely the ones that can be written (modulo a projectivity in $\P^2$) as $\ell_0^{d-2}\ell_1\ell_2$.
\begin{defn} {\rm Let $d\geq 3$; we consider the variety $QQ(V_d) = Im\Phi$, where $\Phi : (\P^2)^*\times (\P^2)^*\times(\P^2)^* \rightarrow \tau_2(V_d)$, with  $\Phi(\ell_0,\ell_1,\ell_2) =F$ and $F\in \P^{N_{d}}$ is the point parameterizing $\ell_0^{d-2}\ell_1\ell_2$; hence $QQ(V_d)$ is a 6-dimensional variety.}
\end{defn}
Note that for $d=2$, $QQ(V_2)=\P^{N_2}=\P^5$.
\dd
{\bf Notation} Let $F$ be a form of degree $d$ on $\P^2$ parameterized by a point $T \in \P^{N_d}$; the symmetric tensor rank $srk(T)$ of $T$ is the $V_d$-rank of $T$ (here $T$ is viewed as a symmetric tensor). 
\begin{prop}\label{MonomialForms} Let $d\geq 3$; we have that if $T\in QQ(V_d)$, $srk(T)\in \{1,d-1,d,2d-2\}$ for $d\neq 3$, while for  $T\in QQ(V_3)$, $srk(T)\in \{1,3,4\}$. In both cases the generic $T$ in $QQ(V_d)$ has $srk(T)=2d-2$. Moreover, $\forall P\in \P^2$:
\begin{enumerate}[leftmargin=*,label=\roman*.]
\item $QQ(V_d)\cap \tau_{2,\nu_d(P)}(V_d) \cong \sigma_2(V_2)=\tau_1(V_2)$;
\item For any $2-$square $Q_P\in \P^2$, $QQ(V_d) \cap \langle\nu_d(Q_P)\rangle \cong \tau_{1,\nu_d(P)}(V_d)\cup \mathcal{Q}_{Q_P}$, where $\mathcal{Q}_{Q_P}\subset \langle\nu_d(Q_P)\rangle \cong \P^3$ is a smooth quadric and we have $\tau_{1,\nu_d(P)}(V_d)=\tau_{1,\nu_d(P)}(\mathcal{Q}_{Q_P})$, i.e. $V_d$ and $\mathcal{Q}_{Q_P}$ have the same tangent plane at $\nu_d(P)$.
\end{enumerate} 
\end{prop}
\proof Let us consider the degenerate cases first: let $T\in \P^{N_d}$ be the point which parameterize a form in $C[x_0,x_1,x_2]_d$ which can be written as $\ell_0^{d-2}\ell_1\ell_2$; if the form is of type $\ell_0^d$, then $T$, of srk =1, belongs to $V_d$; if the form is of type $\ell_0^{d-1}\ell_1$, then $ T \in \tau_1(V_d)$, and more precisely, since it can be written in two variables, there is a rational normal curve $C_d \subset V_d$ such that $T\in \tau_1(C_d)$. Eventually, if the form is of type  $\ell_0^{d-2}\ell_1^2$, then $T\in \tau_2(C_d)$.  It is known (e.g. see \cite{BGI}, Remark 24 or \cite{CCG}, Prop. 3.1 ) that $srk(\ell_0^{d-1}\ell_1)=max\{2,d\}=d$, while $srk(\ell_0^{d-2}\ell_1^2)=max\{3,d-1\}=d-1$, unless $d=3$, when it is =3.  Of course all these $T$'s of the degenerate kind constitute a closed subset $D$ of $QQ(V_d)$. \\
When the form is of type $\ell_0^{d-2}\ell_1\ell_2$, with $\ell_i\neq \ell_j$ $i\neq j \in \{0,1,2\}$, then $srk(T) = 2d-2$ (e.g. see \cite{CCG}). Thus, the part regarding the symmetric rank is proved.
\dd
If we fix the factor $\ell_0^{d-2}$, i.e. we fix the corresponding point on $V_d$ and we consider the osculating space $\tau_{2,\nu_d(P)}(V_d)\cong \P^5$, we have that the points of type $\ell_0^{d-2}\ell_1^2$ are the image under $\Phi$ of the points $(\ell_0,\ell_1,\ell_1)$, i.e. of $\{\ell_0\}\times \Delta$, where $\Delta\cong \P^2$ is the diagonal of $(\P^2)^*\times (\P^2)^*$; hence they form a subvariety $V$ of  $\tau_{2,\nu_d(P)}(V_d)$ isomorphic to $V_2$, the Veronese surface in $\P^5$. We have two ways to check that $QQ(V_d)\cap \tau_{2,\nu_d(P)}(V_d) = \sigma_2(V)$. First, if we consider two points on $V$ (not on $V_d$) parameterizing $\ell_0^{d-2}m_1^2$ and $\ell_0^{d-2}m_2^2$, with $m_1\neq m_2$, we have that the line joining them parameterizes all the forms that can be written as $$\ell_0^{d-2}(\alpha^2m_1^2-\beta^2m_2^2) = \ell_0^{d-2}(\alpha m_1-\beta m_2)(\alpha ml_1+\beta m_2).$$  Since any two lines in a pencil can be projectively transformed in other two lines of the pencil, any form $\ell_0^{d-2}\ell'_1\ell'_2$, with $\ell'_1\neq \ell'_2$ can be projectively transformed into one of the form 
$$\ell_0^{d-2}(\alpha m_1-\beta m_2)(\alpha m_1+\beta m_2),$$
 we get that $QQ(V_d)\cap \tau_{2,\nu_d(P)}(V_d) = \sigma_2(V)$.  Otherwise, and more simply, it suffices to consider that \linebreak $\sigma_2(V_2)=\tau_1(V_2)$, hence $\tau_1(V)$ parameterizes the forms of type $\ell_0^{d-2}\ell_1\ell_2$.  So, $i)$ is proved.
\dd
Now, if we consider a $2-$square $Q_P\in \P^2$, with $I_{Q_P}=(\ell_1^2,\ell_2^2)$, and $\langle\nu_d(Q_P)\rangle \cong  \P^3$, we already know that it contains the plane $\tau_{1,\nu_d(P)}(V_d)\subset D\cap QQ(V_d)$; moreover, we have seen in Prop. \ref{Formequadrellate} that the forms in$\langle\nu_d(Q_P)\rangle$ can be written as $\ell_0^{d-2}(a_0\ell_0^2+a_1\ell_0\ell_1+a_2\ell_0\ell_2+a_3\ell_1\ell_2)$. If we consider its subvariety $\mathcal{Q}_{Q_P}$ given by $a_0a_3-a_1a_2=0$ (a smooth quadric), we have to show that its points parameterize forms in $QQ(V_d)$. Actually, if $a_0=0$, either $a_1=0$ or $a_2=0$, say $a_1=0$, then the form becomes $\ell_0^{d-2}\ell_2(a_2\ell_0+a_3\ell_1)\in QQ(V_d)$; if $a_0\neq 0$, instead, we have that our form becomes $a_0\ell_0^{d-2}(a_0\ell_0+a_1\ell_1)_1(a_0\ell_0+a_2\ell_2)_2$; hence it is in $QQ(V_d)$ again. Note that if $a_3=0$ then the form is $\ell_0^{d-1}(a_0\ell_0+a_1\ell_1+a_2\ell_2)\in \tau_{1,\nu_d(P)}(V_d)$, where either $a_1=0$ or $a_2=0$. Hence $\tau_{1,\nu_d(P)}(V_d)\cap \mathcal{Q}_{Q_P}$ is given by two lines and so the tangent plane to $V_d$ is also tangent to $\mathcal{Q}_{Q_P}$.
\prfend
Let us note that if we knew the equations defining $QQ(V_d)$, we would be able to check if a given form $F\in \C[x_0,x_1,x_2]_d$ can be written as a monomial $x_0^{d-2}x_1x_2$ modulo a linear change of coordinats. Hence it would be interesting to solve the following
\begin{prob} Find equations defining (even just set-theoretically) the variety $QQ(V_d)$.
\end{prob}
\begin{remark} {\rm The above problem could have interest for applications, since we have also that a symmetric tensor in $\P^{N_d}$ describes what in quantum information theory is called a $d-qutrits$ symmetric state, which is not entangled if it is on $V_d$, while its tensor rank can be considered a possible measure of the tensor's entanglement (e.g. see  \cite{BFZ}).}
\end{remark}
\section{$2$-squares on Segre-Veronese surfaces}
\subsection{The first case: $q_2(V_{2,2})$}
We will consider $2$-squares also in $\P^1\times \P^1$ and their Segre-Veronese embeddings $\nu_{d,d}$, given by the forms of bidegree $(d,d)$, $d\geq 2$, in $R=\C[s_0,s_1;t_0,t_1]$. 
 We are not to consider all the possible 2-$squares$ in $\P^1\times \P^1$, we consider only the simplest ones, i.e. the ones whose ideal is generated in bidegrees $(2,0)$ and $(0,2)$. Given bihomogeneous coordinates $[s_0,s_1;t_0,t_1]$ in  $\P^1\times \P^1$, consider the affine chart $U_{1,1}=\{s_1\neq 0, t_1\neq 0\}$ of $\P^1\times \P^1$, with coordinates $(x,y)$, where $x =\frac{s_0}{s_1} ,y = \frac{t_0}{t_1}$; for any point $P=[a,b;c,d]\in U_{1,1}$ we can consider the 0-dimensional scheme $Q_P$ defined by the bihomogeneous ideal $(\ell_s^2,\ell_t^2)\subset \C[s_0,s_1;t_0,t_1]$, where $\ell_s = bs_0-as_1$ and $\ell_t=dt_0-ct_1$. Since $bd\neq 0$ we can look at $Q_P$ in the affine chart of coordinates $(x,y) = (s_0/s_1,t_0/t_1)$ and there $I_{Q_P} = ((x-a/b)^2, (y-c/d)^2)$.  So, $Q_P$ is a 2-square in the affine chart $U_{1,1}\cong \A^2$. This will be the kind of $2-$squares that we are going to consider in multi-projective enviroment. 
 \dd
Now we consider the Segre-Veronese embedding $sv_{2,2}: \P^1\times \P^1 \rightarrow \P^8$ given by $\mathcal{O}_{\P^1\times \P^1}(2,2)$, i.e. 
$$sv_{2,2}(s_0,s_1;t_0,t_1) = (s_0^2t_0^2,2s_0^2t_0t_1, s_0^2t_1^2, 2s_0s_1t_0^2,4s_0s_1t_0t_1, 2s_0s_1t_1^2,s_1^2t_0^2,2s_1^2t_0t_1, s_1^2t_1^2).$$
The image  $V_{2,2} = sv_{2,2}(\P^1\times \P^1)\subset \P^8$ is the (2,2)-Segre-Veronese surface of $\P^1\times \P^1$.  Using coordinates $[z_{i_0,i_1,j_0,j_1}]$, $i_0,i_1,j_0,j_1 \in \{0,1\}$, $i_0+i_1=2$, $j_0+j_1=2$, $i_0\leq i_1$, $j_0\leq j_1$, in this $\P^8$, the parametric equations of $V_{2,2}$ are  $z_{i_0,i_1,j_0,j_1}=(s_0^{i_0}s_1^{i_1}t_0^{j_0}t_1^{j_1})$. We have that $V_{2,2}$ is a particular Del Pezzo surface, of degree 8 in $\P^8$ (e.g. see \cite{CGG2}); its ideal is generated by the $2\times 2$ minors of the matrix 
$$\left( \begin{array} {cccc} z_{0000}&z_{0001} &z_{0100}&z_{0101}\cr 
 z_{0001}&z_{0011} &z_{0101}&z_{0111} \cr 
 z_{0100}&z_{0101} &z_{1100}&z_{1101}\cr z_{0101}&z_{0111} &z_{1101}&z_{1111} \end{array}\right)$$
Moreover the $3\times 3$-minors of the matrix above generate its secant variety, $\sigma_2(V_{2,2})$, which has the expected dimension =5, while its determinant defines $\sigma_3(V_{2,2})$, which is defective, since its expected dimension was 8 (e.g. see again \cite{CGG2}).	
\dd
$V_{2,2}$ can be seen as given first by the 2-ple Veronese embedding of both the $\P^1$-factors into $\P^2$, followed by the Segre embedding  $\nu_{1,1}:\P^2\times \P^2\rightarrow \P^8$. Moreover, this space $\P^8$  parameterizes $2^4$-tensors which are symmetric on the first two indices and on the second two; namely we can see these partially symmetric tensors as the subspace ($\cong \P^8$) of the space of general $2^4$-tensors ($\cong \P^{15}$) with entries $x_{i_0,i_1,u_0,u_1,j_0,j_1,v_0,v_1}$  defined by the equations  $x_{i_0,i_1,u_0,u_1,j_0,j_1,v_0,v_1} =x_{\sigma((i_0,i_1),(u_0,u_1)),\tau((j_0,j_1),(v_0,v_1)}$, for all $\sigma,\tau \in \mathfrak{S}_2$. We can view the partially symmetric $2^4$-tensors in the following figure. Note that the four $2\times 2$ faces moving from left to right are symmetric and so are the four $2\times 2$ faces joining the ``big cube" to the ``small one" in the direction ``perpendicular to the paper"  
\begin{figure}[H]
		\centering
		\includegraphics[scale=0.35]{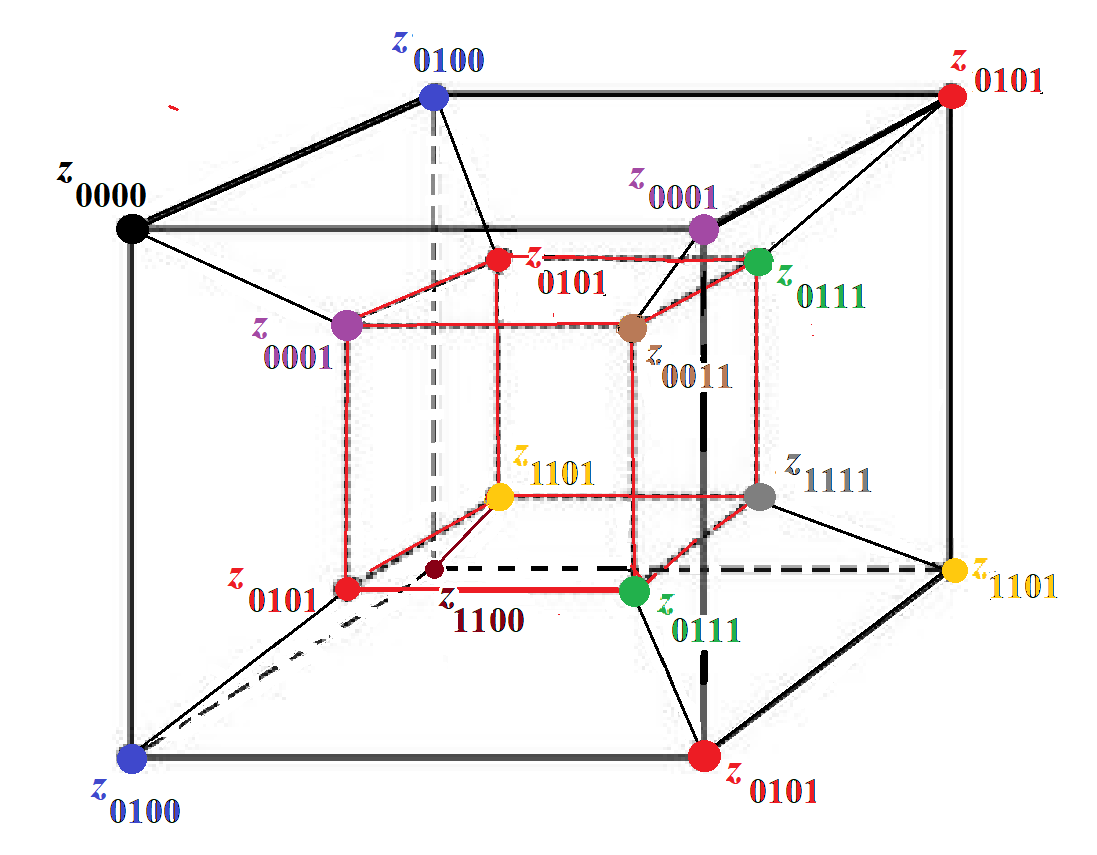}
		\caption {Partially symmetric $2^4$-tensors.}
	\end{figure}
\noindent The generators of $I_{V_{2,2}}$ can also be viewed as the $(2\times 2)-$minors of the array above.
\dd
{\bf Notation}  The partially symmetric tensor rank $psrk(T)$ of a point $T\in \P^{8}$ is the $V_{2,2}-rank$ of $T$. 
\bigskip
Now let us go back to the $2$-square schemes we defined in $\P^1\times \P^1$; we considered a scheme $Q_P$ for every point $P\in \P^1\times \P^1$; under the Segre-Veronese embedding whose image is $V_{2,2}\subset \P^8$, each of these schemes, which we will denote $Q'_P$, is such that $\langle Q'_P\rangle \cong \P^3$. Here we will consider the variety  
$$q_2(V_{2,2}) = \overline{\bigcup_{P\in \P^1\times \P^1} \langle Q'_P\rangle}$$ 
In order to visualize what $\langle Q'_P\rangle$ is, recall that we defined $I_{Q_P} = (\ell_s^2,\ell_t^2)\subset R=\C[s_0,s_1;t_0,t_1]$; let \linebreak $m_s = as_0+bs_1\in R_{1,0}$ and $m_t = ct_0+dt_1\in R_{0,1}$, those are the elements of bidegree (1,0) and (0,1) such that $m_s\perp \ell_s$, $m_t\perp \ell_t$, with respect to the apolar action of the ring on itself via derivations (note that we can have $m_s=\ell_s$ or $m_t=\ell_t$ for isotropic forms).  We have that 
$$(I_{Q_P})_{(2,2)} = \ell_s^2R_{0,2} + \ell_t^2R_{2,0},$$
so that
$$(I^\perp_{Q_P})_{(2,2)} = m_sm_tR_{1,1} \cong \ \langle Q'_P\rangle.$$
Let us consider for example $P = [0,1;0,1]$, then $I_{Q_P} = (s_0^2,t_0^2)$ and  
$$(I^\perp_{Q_P})_{(2,2)} = \langle s_0s_1t_0t_1,s_0s_1t_1^2,s_1^2t_0t_1,s_1^2t_1^2\rangle,$$ 
hence the tensors in any $\langle Q'_P\rangle$ can be written, modulo a bilinear change of coordinates in $\P^1\times \P^1$, as described in the next figure:
\begin{figure}[H]
		\centering
		\includegraphics[scale=0.35]{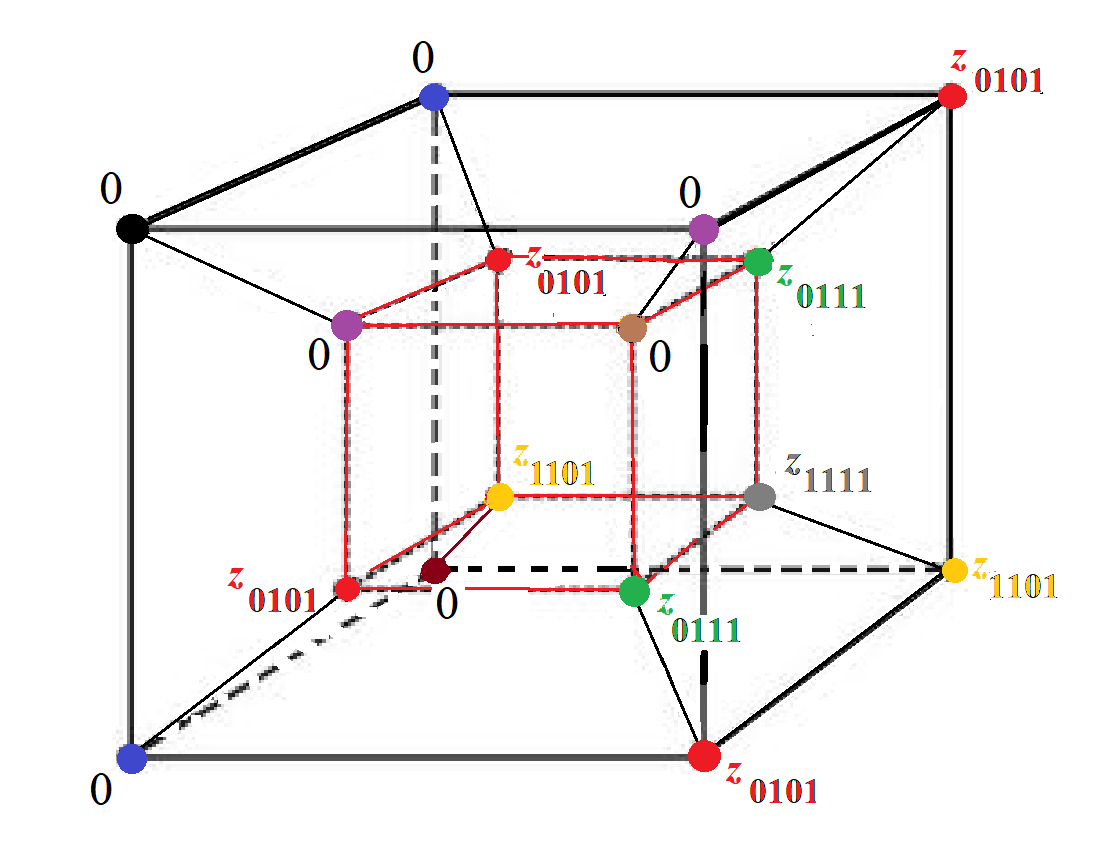}
		\caption {Tensors in $\langle Q'_P\rangle$}
	\end{figure}
\noindent Hence $\langle Q'_P\rangle$ is defined by equations $z_{i_0i_1j_0j_i}=0$, for all $z_{i_0i_1j_0j_i}$ where either both $i$'s or both $j$'s are 0.  Note that since $\langle Q'_P\rangle = \P(m_sm_tR_{1,1})$, while $\tau_{1,\nu_{22}(P)} = \P\langle m_sm_t(m_sR_{0,1}+m_tR_{1,0})\rangle$, we have that \linebreak $\tau_1(V_{2,2})\subset q_2(V_{2,2})$.  

\medskip
Now we want to consider the secant variety of $q_2(V_{2,2})$; namely we want to prove the following.

\begin{prop} We have that $\dim q_2(V_{2,2}) =5$ and $\sigma_2(q_2(V_{2,2})) = \P^8$, as expected. Hence the generic partially symmetric tensor in $\P^8$ can be written as the sum of two partially symmetric tensors which depend only on four parameters each (and can be written, not at the same time, as in Fig. 2).
\end{prop}
\proof To give a point in $\langle Q'_P\rangle = \P^3 \subset \P^8$ amounts to choosing a form $m_sm_tn_{s,t}$, with $m_s \in R_{1,0}$, $m_t\in R_{0,1}$, $n_{s,t}\in R_{1,1}$.  Hence, in order to find the tangent space to $q_2(V_{2,2})$ at the point $m_sm_tn_{s,t}$, we have to consider another generic point $u_su_tv_{s,t} \in q_2(V_{2,2}) $, and then compute (e.g. see \cite{CGG}, \cite{CGG1}):
$$ \lim _{\lambda \rightarrow 0}	\frac{d}{d\lambda}\left[ (m_s+\lambda u_s)(m_t+\lambda u_t)(n_{s,t}+\lambda v_{s,t})\right] =$$ $$ = u_sm_tn_{s,t} + m_su_tn_{s,t} +m_sm_tn_{s,t} \subset R_{2,2}.$$

As $u_s , u_t, v_{s,t}$ vary, we get that the affine cone on the tangent space that we considered is:
$$ W = m_sm_tR_{1,1} + m_sn_{s,t}R_{0,1} + m_tn_{s,t}R_{1,0} . $$ 
We have dim $W$ = 6 (affine), as expected.  If we consider the ideal $I = (m_sm_t, m_sn_{s,t}, m_tn_{s,t})$, we have that  $W = I_{2,2}$, and $I$ is the ideal of three points $P_1,P_2,P_3 \in \P^1\times \P^1$; since they can be defined by $(m_s,n_{s,t})$, $(m_t,n_{s,t})$ and $(m_s,m_{t})$, respectively, the three of them are not contained in a fiber, but $P_1,P_3$ and $P_2,P_3$ are (see Fig. 3a). 
\begin{figure}[H]\label{3Punti}
		\centering
		\includegraphics[scale=0.65]{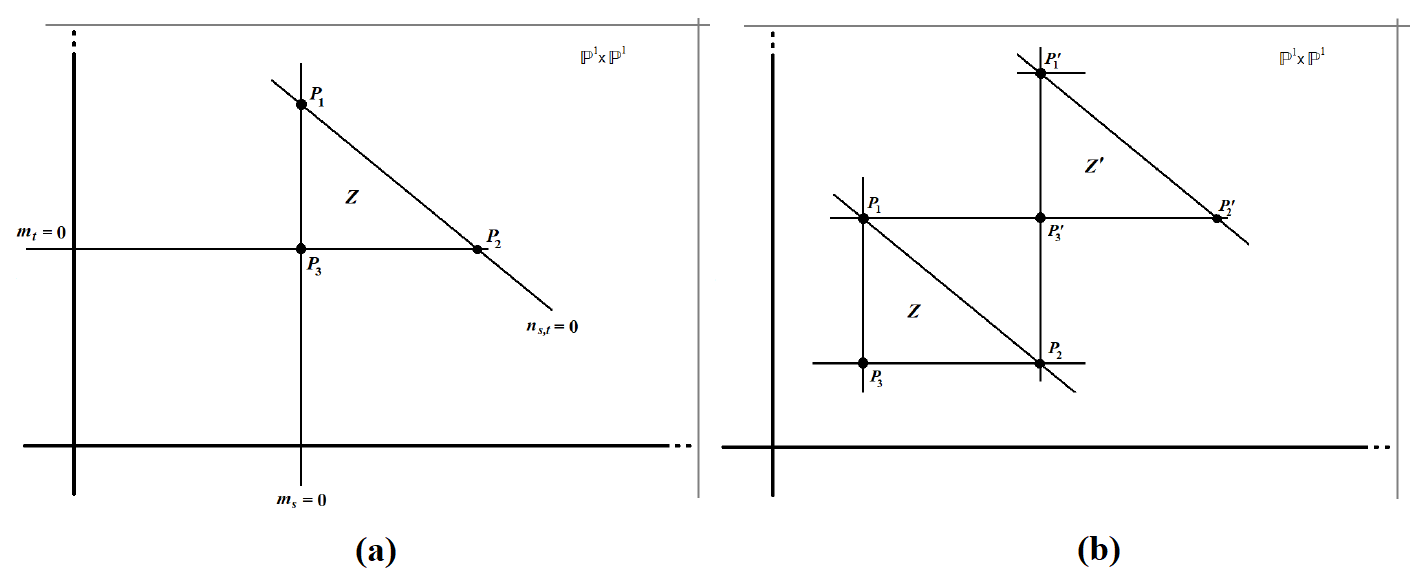}
		\caption {{\bf (a)} $P_1,P_2,P_3$ in $\P^1\times \P^1$; {\bf (b)} $Z\cup Z'$ specialized. }
	\end{figure} 

We just have to consider two tangent spaces to  $q_2(V_{2,2})$ at two generic points of $q_2(V_{2,2})$ (say given by $m_sm_tn_{s,t}$ and $u_su_tv_{s,t}$); if their affine cones are $W$ and $W'$, the (affine) space $W+W'$ will give the tangent cone at a generic point of  $\sigma_2(q_2(V_{2,2}))$, by Terracini's Lemma. 
Since $W=(I_Z)_{2,2}$ and $W'=(I_{Z'})_{2,2}$, where $Z$ and $Z'$ are given by two configuration of three simple points as in Fig. 3.b, $W\cap W'$ is the $(2,2)$ part of the ideal of $Z\cup Z'$.

\medskip
\noindent {\bf Claim}  {\it Let $Z,Z' \subset \P^1\times \P^1$ be two schemes of three points,  positioned as in Fig. 3b. Then they impose independent conditions to forms of bidegree $(2,2)$, i.e. $\dim (I_{Z'\cup Z})_{2,2} = 3$.}

\medskip
\noindent{\bf Proof of Claim} We can specialize $Z'$ (see Fig. 6b) so that the line $u_s=0$ contains $P'_3$, $P'_1$ and also a point (say $P_2$) of $Z$; this forces the forms in $(I_{Z'\cup Z})_{2,2}$ to be of type $u_sF$, with $F\in R_{1,2}$; now we have also that $F$, since $\{F=0\}$ contains $P_1$ and $P'_2$, has to be of type $F=u_tG$, with $G\in R_{1,1}$ and $P_3\in \{G=0\}$. Hence $\dim  (I_{Z'\cup Z})_{2,2} = \dim (I_{P_3})_{1,1} = 3$.\prfend
\dd
By the claim we get $\dim (W\cap W') = \dim (I_{Z\cup Z'})_{2,2} = 3$, so 
$$\dim(W+W') = \dim W + \dim W' - \dim (W\cap W') = 6 + 6 - 3 = 9$$
 and $\sigma_2(q_2(V_{2,2})) = \P^8$. \prfend
\subsection{The varieties $V_{d,d}$ and their $q_2(V_{d,d})$, $d\geq 3$.}
Now we want to generalize what we did to the case of Segre-Veronese varieties $$V_{d,d}=\nu_{d,d}(\P^1\times\P^1)\subset \P^{(d+1)^2-1},$$
 with $d\geq 3$. Also in this case $V_{d,d}$ parameterizes partially symmetric $2^{2d}-$tensors; we can define $q_2(V_{d,d})$ exactly as before and we have:
\begin{prop} For all $d\geq 3$, dim $q_2(V_{d,d}) = 5$ and dim $\sigma_2(q_2(V_{d,d})) = 11$, as expected.
\end{prop}
\proof We consider the variety $V_{d,d}\subset \P^{(d+1)^2-1}$, given by the Segre-Veronese embedding $\nu_{d,d}$, associated to the very ample line bundle $\mathcal{O}_{\P^1\times\P^1}(d,d)$; here $\P^{(d+1)^2-1} \subset \P^{2^{2d}-1}$, parameterizes the $2^{2d}$-tensors which are partially symmetric in the sense that each $z_{i_{1},\ldots ,i_d,i_{d+1},\ldots ,i_{2d}} = z_{\sigma(i_1,\ldots,i_d),\sigma'(i_{d+1},\ldots,i_{2d})}$, for all $\sigma , \sigma'\in \mathfrak{S}_d$. The variety $V_{d,d}$ parameterizes the partially symmetric $2^{2d}-$tensors $T$ with $psrk T =1$; as in the case $d=2$, the partially symmetric tensor rank $psrk(T)$ of a point $T\in P^{(d+1)^2-1}$ is the $V_{d,d}$-rank of $T$.
\dd
We want to consider the variety $q_2(V_{d,d})$, given by all the $2$-squares $Q_P$ defined by ideals $(\ell_s^2,\ell_t^2)$ and their images $Q_P' \subset V_{d,d}$.  Here $\langle Q'_P\rangle$ will correspond to the linear space given by 
$$(I^{\perp}_{Q_P})_{d,d} =  m_s^{d-1}m_t^{d-1}R_{1,1} = (m_s^{d-1}m_t^{d-1})_{d,d}.$$
Hence, in order to find the tangent space to $q_2(V_{d,d})$ at  
the point corresponding to $m_s^{d-1}m_t^{d-1}n_{s,t}$, we have to consider another generic point $u_s^{d-1}u_t^{d-1}v_{s,t} \in q_2(V_{d,d}) $, and then to compute:
$$ \lim _{\lambda \rightarrow 0}	\frac{d}{d\lambda}\left[ (m_s+\lambda u_s)^{d-1}(m_t+\lambda u_t)^{d-1}(n_{s,t}+\lambda v_{s,t})\right] =$$ 
$$ =(d-1)m_s^{d-2}u_sm_t^{d-1}n_{s,t} + (d-1)m_t^{d-2}u_tm_s^{d-1}n_{s,t}+ m_s^{d-1}m_t^{d-1}v_{s,t}.$$
This, as $u_s,u_t,v_{s,t}$ vary, gives the space
$$W = (m_s^{d-1}m_t^{d-1}, m_s^{d-1}m_t^{d-2}n_{s,t}, m_s^{d-2}m_t^{d-1}n_{s,t})_{d,d} \subset (I_Z)_{d,d}$$
where $Z$ is the scheme which is made of the two lines $\{m_s=0\}$, $\{m_t=0\}$, both with multiplicity $d-2$, plus two $(d-1)$-jets on $\{n_{s,t}=0\}$, supported at the points given by $(m_s,n_{s,t})$ and $(n_{s,t},m_t)$ (see Fig.4). We have that $W$ has (vector) $dim W = 6$, as expected, while $dim(I_Z)_{d,d} = 7$ (all forms in $(I_Z)_{d,d}$ are of type $m_s^{d-2}m_t^{d-2}F$, where $F$ is a $(2,2)-$form passing through the two points where the jets are, hence $dim(I_Z)_{d,d} = 9-2=7$). 
\begin{figure}[H]\label{(d-2)lines(d-1)jets}
		\centering
		\includegraphics[scale=0.6]{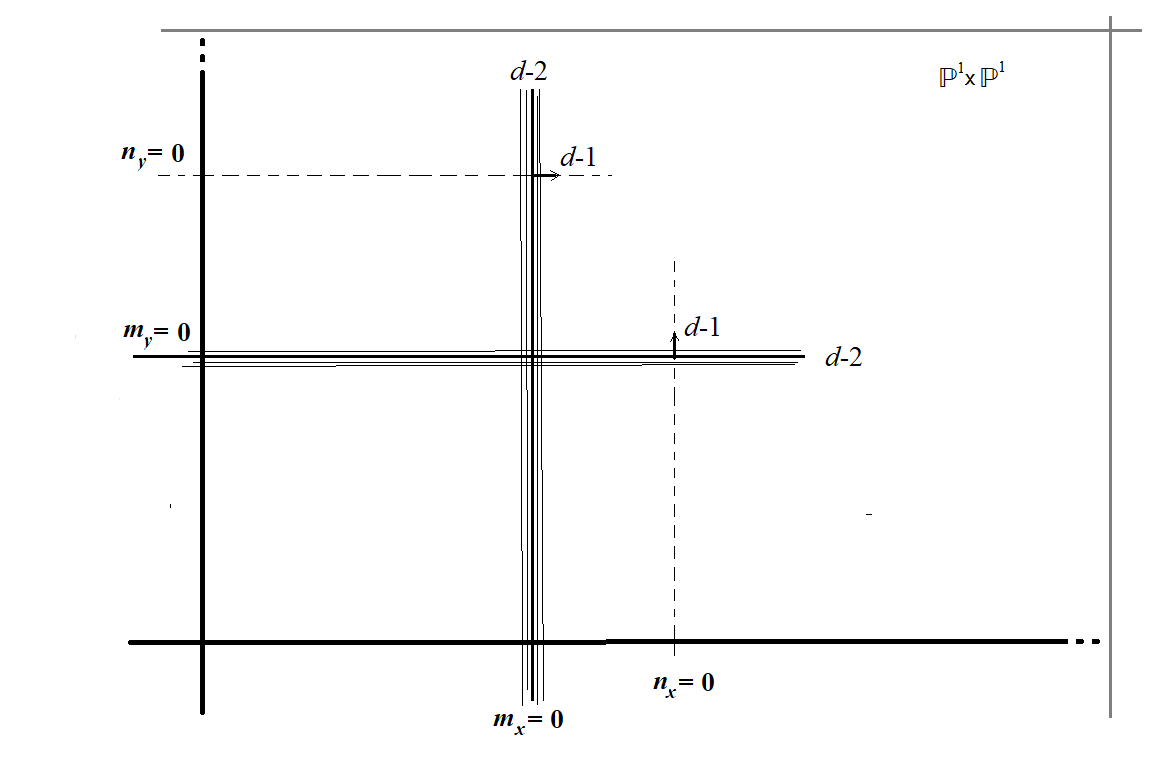}
		\caption {The scheme $Z$: two $(d-2)-$ple lines plus two $(d-1)$-jets on a third line.}
	\end{figure}
\noindent By Terracini's Lemma, we just have to consider two tangent spaces to $q_2(V_{d,d})$ at two generic points of $q_2(V_{d,d})$; if their affine cones are $W$ and $W'$, the (vector) space $W+W'$ will be the tangent cone to a generic point of  $\sigma_2(q_2(V_{3,3}))$. 
Since $W\subset (I_Z)_{d,d}$ and $W'\subset (I_{Z'})_{d,d}$, where $Z$ and $Z'$ are made as in Fig. 4 (two $(d-2)-$lines with two $(d-1)-$jets), $W\cap W' \subset (I_Z)_{d,d}\cap (I_{Z'})_{d,d} = (I_{Z\cup Z'})_{d,d}$.  We want to check that $(I_{Z\cup Z'})_{d,d} =\{0\}$.
\begin{itemize}[leftmargin=*]
	\item Case $d=3$: since the forms in $(I_{Z\cup Z'})_{3,3}$ should contain the factor $m_sm_tm'_sm'_t\in R_{2,2}$, and no form in $R_{1,1}$ passes through the four points which are the support of the four jets, we get $(I_{Z\cup Z'})_{3,3}= \{0\}$;
	\item Case $d\geq 4$: forms in $(I_{Z\cup Z’})_{d,d}$  should contain the factor $m_s^{d-2}m_t^{d-2}m_s'^{d-2}m_t'^{d-2}$, which is impossible for $d\geq 5$, while for $d=4$ there would be only this form, which does not contain the four $d-1$ jets at the four points $m_s \cap n_{s,t}$ , $m_t \cap n_{s,t}$  , $m’_s \cap n'_{s,t}$ , $m’_t \cap n'_{s,t}$, hence also in this case $(I_{Z\cup Z'})_{3,3}= \{0\}$, which in turn implies that $W\cap W' = \{0\}$.  
\end{itemize}
Now, from Grassmann equality, we get $\dim (W+W') = \dim W + \dim W' - \dim (W\cap W') = 6+6-0$, hence $\dim (W+W') =12$ (as vector space) and $\dim\sigma_2(q_2(V_{3,3})) = 11$. \prfend
\subsection{The varieties $qq_2(V_{d,d})$}
As we have done in the Veronese case, here too we are going to define a subvariety of $q_2 (V_{d,d})$, parameterizing a special kind of partially symmetric tensors which could also have some interest in relation to Quantum Entanglement. 
If we consider the Hilbert space of a composite quantum system, then this is the tensor product of the Hilbert spaces of the constituent systems, and tensor rank is a natural measure of the entanglement for the corresponding quantum states. The Hilbert space of a $k$-body system is obtained as the tensor product of $k$ copies of the single particle Hilbert space ${\mathcal H}_1$. In the case of indistinguishable bosonic particles, the totally symmetric states under particle exchange are physically relevant, which amounts to restricting the attention to the subspace $H_s = S^N({\mathcal H}_1)\subset  {\mathcal H}_1^{\otimes N}$ of symmetric tensors.
\dd
In case we have $k\geq 2$ different species of indistinguishable bosonic particles, the relevant Hilbert space is  $S^{N_1}({\mathcal H}_1)\otimes  S^{N_2}({\mathcal H}_2)\otimes \ldots \otimes S^{N_k}({\mathcal H}_2)$.  Of particular interest, in physics literature, are the so-called $W$-states, i.e. quantum entangled states that can be expressed in Dirac notation as:
$${\displaystyle |\mathrm {W} \rangle ={\frac {1}{\sqrt {n}}}(|100...0\rangle +|010...0\rangle +...+|00...01\rangle ).}$$ 
Which, if ${\mathcal H}_1 = \C^2$ with coordinates $(x,y)$, can be written as:
$$W_d = y\otimes x\otimes \ldots \otimes x +  x\otimes y\otimes x\otimes \ldots \otimes x + x\otimes x\otimes \ldots \otimes y.$$
When treating with bosonic particles, hence with symmetric tensors, $W_d$ can be represented simply as a monomial $x^{d-1}y$, hence in the study of entanglement of $k$ different $d$-body systems (made of different species of indistinguishable bosonic particles, like photons), we can consider the product of $k$ $W_d$ states: $W_d \otimes \ldots \otimes W_d$, where each $W_d \cong x^{d-1}y \in S^d \C^2  \subset  (\C^2)^{\otimes d}$ (e.g. see \cite{BBCG}). In  \cite{BBCG}, Lemma 2.1, it is proved that $W_d \otimes \ldots \otimes W_d \in \langle \nu_{d,\ldots,d}(Z)\rangle$, where $Z$ is a $2-$hypercube in $\P^1\times \ldots \times \P^1$, and this is expressed by saying that the cactus rank of $W_d \otimes W_d$ is $2^k$ and it is realized by $Z$ (the cactus rank of a tensor $T\in \P^{(d+1)^k-1}$ is the minimum lenght of a $0-$dimensional scheme $X\subset V_{d,\ldots ,d}$ such that $T\in \langle X\rangle$). We can improve a bit that lemma in this setting (in \cite{BBCG} also the case of  $W_{d_1} \otimes \ldots \otimes W_{d_k}$ is considered, with different $d_i$'s).
\begin{cor} Let $T\in P^{(d+1)^k-1}$ parameterize $W_d \otimes \ldots \otimes W_d$; then the smoothable rank of $T$ is $smrk = 2^k$. 
\end{cor}
\proof The only difference between smoothable rank and cactus rank is that $smrk\ T = r$ if and only if there is a smoothable  $0-$dimensional scheme $X\subset V_{d,\ldots ,d}$ such that $T\in \langle X\rangle$ and $\ell(X) =r$. Since, by Proposition \ref{smoothable}, any $2-$hypercube is smoothable, the statement is an immediate consequence of this and of Lemma 2.1 of \cite{BBCG}, since for any tensor the smoothable rank is greater or equal than the cactus rank. \prfend
Let us consider the case  $k=2$;  $q_2(V_{d,d})$ parameterizes partially symmetric tensors in the spaces:
$$<\nu_{d,d}(Q_P)> = (m_s^{d-1}m_t^{d-1})_{d,d}$$ 
Hence for all tensors of type  $W_d \otimes W_d = m_s^{d-1}a_sm_t^{d-1}b_t$, $m_s,a_s \in R_{1,0}$, $m_t,b_t \in R_{0,1}$, we have  \linebreak $W_d \otimes W_d\in <\nu_{d,d}(Q_P)>$, for some $P\in \P^1\times \P^1$.  
More specifically, let us consider the subvariety  which parameterizes exactly the tensors of type $W_d \otimes W_d$:
\begin{defn} {\rm The} variety  {\rm $qq_2(V_{d,d}) \subset  q_2(V_{d,d})$ is the image of the map:
$$ \P^1 \times \P^1 \times \P^1 \times \P^1 = \P(\C[s_0,s_1]_1)\times \P(\C[s_0,s_1]_1)\times\P(\C[t_0,t_1]_1)\times\P(\C[t_0,t_1]_1)\rightarrow   q_2(V_{d,d}) \subset \P^8,$$ 
with  $(m_s, n_s, m_t, n_t ) \rightarrow   m_s^{d-1}n_sm_t^{d-1}n_t$.} 
\end{defn}
We have that $qq_2(V_{d,d})$  has dimension 4, and, via a multi-linear change of coordinates, every form parameterized by a point in $qq_2(V_{d,d})$ can be written as a monomial $s_0^{d-1}s_1t_0^{d-1}t_1$ (for results on the various ranks of such monomials see \cite{CF},\cite{Ga} and \cite{BBCG}). We have the following proposition.
\begin{prop}\label{qq} $qq_2(V_{d,d})$ is such that:
\begin{enumerate}[label=\roman*.,leftmargin=*]
\item $\forall P \in \P^1\times \P^1$, $qq_2(V_{d,d})\cap \langle Q'_P\rangle \cong \mathcal{Q}_P$, where $\mathcal{Q}_P$ is a smooth quadric in $\langle Q'_P\rangle \cong \P^3$.
\item $\forall P \in \P^1\times \P^1$, we have $\tau_{1,\nu_{dd}(P)}(\mathcal{Q}_P) = \tau_{1,\nu_{dd}(P)}(V_{dd})$.
\item$Sing(qq_2(V_{d,d}))$ is the locus of forms of type $m_s^{d}m_t^{d-1}n_t$ or $m_s^{d-1}n_sm_t^{d}$.
\end{enumerate}
\end{prop}
\proof As we have seen, if $\nu_{d,d}(P) = m_s^dm_t^d$, then we can write   
$$\langle Q'_P\rangle = (m_s^{d-1}m_t^{d-1})_{d,d} =  m_s^{d-1}m_t^{d-1}R_{1,1},$$
where $n_s,n_t$ are such that $R_{1,0}=\langle m_s,n_s\rangle$, $R_{0,1}=\langle m_t,n_t\rangle$ and $R_{1,1}= \langle m_sm_t, m_sn_t, n_sm_t, n_sn_t\rangle$. Hence any points in $\langle Q'_P\rangle$ corresponds (modulo constants) to a form of type $am_sm_t+ bm_sn_t + cn_sm_t + dn_sn_t$, i.e. we can view $\langle Q'_P\rangle \cong \P(R_{1,1})$, and in it the forms of type $a_sa_t$, for $a_s\in R_{1,0}$, $a_t\in R_{0,1}$, are precisely parametrized by a quadric $\mathcal{Q}_P$ which is the Segre Variety $V_{11}$; $\mathcal{Q}_P$ is given exactly by the forms for which $ad-bc=0$. This proves part  {\it i)}.
\dd
To prove part {\it ii)}, just note that $\tau_{1,\nu_{dd}(P)}(V_{dd})$ is given by the forms in $m_s^{d-1}m_t^{d-1}(m_sR_{0,1}+ m_tR_{1,0})$, hence $\mathcal{Q}_P\cap \tau_{1,\nu_{dd}(P)}(V_{dd})$ is given by the forms of type either $m_s^dm_t^{d-1}(\alpha_t+m_t)$ or $m_s^{d-1}(m_s+\beta_s)m_t^{d}$, which gives two lines in $\tau_{1,\nu_{dd}(P)}(V_{dd})$, hence this is the tangent plane to $\mathcal{Q}_P$ in $\langle Q'_P\rangle $.
\dd
In order to prove  {\it iii)}, let us consider the affine cone $W$ over the tangent space of $qq_2(V_{d,d})$ at one of its points, say the one associated to $m_s^{d-1}n_sm_t^{d-1}n_t$; if we consider another point $u_s^{d-1}v_su_t^{d-1}v_t$, we have to compute:
$$  \lim _{\lambda \rightarrow 0}	\frac{d}{d\lambda}\left[(m_s+\lambda u_s)^{d-1}(n_s+\lambda v_s)(m_t+\lambda u_t)^{d-1}(n_st+\lambda v_t)\right] =$$
$$= \langle (d-1)m_s^{d-2}u_sn_sm_t^{d-1}n_t + m_s^{d-1}v_sm_t^{d-1}n_t + (d-1)m_s^{d-1}n_sm_t^{d-2}u_tn_t+m_s^{d-1}n_sm_t^{d-1}v_t\rangle .$$
Hence, as $u_s,v_s,u_t,v_t$ vary, we get 
$$W \cong \langle m_s^{d-2}m_t^{d-2}\left( m_sm_t(R_{1,0}n_t + n_sR_{0,1}) + n_sn_t(R_{1,0}m_t + m_sR_{0,1})\right)\rangle \subset R_{d,d}. $$
Generically, we have $\dim_\C (R_{1,0}a_t + a_sR_{0,1}) = 3$, since they have $\langle a_sa_t\rangle$ in common, hence $W$ is the sum of two subspaces of (affine) dimension 3, which have $\langle m_s^{d-1}n_sm_t^{d-1}n_t\rangle$ in common, so $\dim_\C W = 5$, as expected. The locus $Sing(qq_2(V_{d,d}))$ is given by the points where $\dim W < 5$, and it is easy to check that this happens for either $m_s=n_s$ or $m_t=n_t$, and this proves  {\it iii)}.
\prfend
There is another way to view the variety $qq_2(V_{d,d})$; consider the embedding $\nu_{dd}(\P^1\times \P^1)$ as the composition:
$$\P^1\times \P^1 \rightarrow \P^d\times \P^d \rightarrow \P^{d^2+2d}$$
where the first arrow is $\nu_d\times \nu_d$ and the second is the Segre embedding $s_{1,1}$. If the image of the first map is $C_d^s\times C_d^t$, where $C_d^s, C_d^t$ are the rational normal curves defined by $\C[s_0,s_1]_d, \C[t_0,t_1]_d$, respectively, we can consider the product of their tangential varieties $\tau(C_d^s)\times  \tau (C_d^t)\subset \P^d\times \P^d$, parameterizing pairs of forms like $(m_s^{d-1}a_s,m_t^{d-1}a_t)$; so $s_{1,1}(\tau(C_d^s)\times  \tau (C_d^t))$ is exactly $qq_2(V_{d,d})$. We know that the singular locus of the tangential surface to a rational normal curve is the rational normal curve itself, hence 
$$Sing(\tau(C_d^s)\times  \tau (C_d^t)) = (\tau(C_d^s)\times  C_d^t) \cup (C_d^s\times  \tau (C_d^t)),$$ in correspondence with what we saw in Prop. \ref{qq}, {\it iii)}.
\dd
Note that for $d=2$, we have:  $qq_2(V_{2,2}) = s_{1,1}(\tau(C^s_2)\times \tau(C^t_2)) = s_{1,1}(\P^2\times \P^2)$, hence $qq_2(V_{2,2})$ is just the Segre Variety $S = s_{1,1}(\P^2\times \P^2) \subset \P^8$, which is well-known to be $2-$defective (it is the variety of $3\times 3$ matrices of rank 2), i.e. dim $\sigma_2(S) = 7$. We want to check that this does not happen for $d\geq 3$, i.e.
\begin{prop}
For $d\geq 3$, $\dim \sigma_2(qq_2(V_{d,d})) = 9$, as expected.
\end{prop}
\proof  By Terracini's Lemma the dimension of the affine tangent cone at a generic point of $\sigma_2(qq_2(V_{d,d}))$ will be $\dim W_1+W_2$, where $W_1$, $W_2$ are the affine tangent cones at two generic points of $qq_2(V_{d,d})$. Thus, in order to prove our statement, we have to show that $\dim(W_1+W_2) = 10$, i.e. since $\dim W_1 = \dim W_2 = 5$, that $W_1\cap W_2 = \{ 0\}$.\\
In the proof of Prop. \ref{qq}, {\it iii)}, we have computed the affine tangent cone $W$ at a generic point of $qq_2(V_{2,2})$, hence if we pick two generic points given by forms: $m_s^{d-1}n_sm_t^{d-1}n_t$ and $u_s^{d-1}v_su_t^{d-1}v_t$, we will have:
$$W_1 = \langle m_s^{d-2}m_t^{d-2}\left( m_sm_tn_tR_{1,0} + m_sm_tn_sR_{0,1} + n_sm_tn_tR_{1,0} + n_sm_sn_tR_{0,1}\right)\rangle \subset ( m_s^{d-2}m_t^{d-2})R_{2,2},$$
$$W_2 = \langle u_s^{d-2}u_t^{d-2}\left( u_su_tv_tR_{1,0} + u_sv_su_tR_{0,1} + v_su_tv_tR_{1,0} + u_sv_sv_tR_{0,1}\right)\rangle \subset ( u_s^{d-2}u_t^{d-2})R_{2,2}.$$
When $d\geq 5$ it is immediate to check that $W_1\cap W_2 = \{ 0\}$, so we are done.
\dd
If $d=4$, then  $( m_s^{d-2}m_t^{d-2})R_{2,2}\cap ( u_s^{d-2}u_t^{d-2})R_{2,2} = \langle m_s^{2}m_t^{2}u_s^{2}u_t^{2}\rangle$ and it is easy to check that $m_s^{2}m_t^{2}u_s^{2}u_t^{2} \notin W_1\cap W_2$ (it sufficies to consider $m_s = s_0, m_t = t_0, n_s = s_1, n_t =t_1$ and note that $s_i^2t_j^2 \notin W_1$, while these monomials will appear in $u_s^{2}u_t^{2}$).
\dd
For $d=3$, let $W_1 = s_0t_0(s_0^2,s_0s_1t_0t_1,s_0s_1t_0^2,s_1^2t_0t_1,s_0s_1t_1^2)$, if there is something not 0 in $W_1\cap W_2$ it should be of the form $s_0t_0u_su_tn_{s,t}$, with $n_{s,t}\in R_{1,1}$, but since $u_s,u_t$ are generic, say  $u_su_t = (as_0+bs_1)(ct_0+bt_1)$, in $u_su_tn_{s,t}$ the monomials $s_i^2t_j^2$ should appear, and this is impossible since they are not in $\langle s_0^2,s_0s_1t_0t_1,s_0s_1t_0^2,s_1^2t_0t_1,s_0s_1t_1^2 \rangle$. Hence  $W_1\cap W_2 = \{ 0\}$ also for $d=3$ and we are done.
\prfend
\bigskip
\noindent {\bf Acknowledgements} All the authors are members of GNSAGA of Indam. We would like to thank J. Buczy\'nski for his several suggestions. The third author wants to thank the organizers of AGATES workshop “Geometry of secants”, Warsaw October 2022 (supported by the Simons Foundation, award number 663281, and by Institute of Mathematics of Polish Academy of Sciences) for the opportunities given by this useful initiative. In particular thanks to  H. Abo, E.Ventura and E. Postinghel for useful talks.

\medskip \b {\bf Funding }Stefano Canino, Alessandro Gimigliano and Monica Idà have been funded by the \lq\lq 0-Dimensional Schemes, Tensor Theory, and Applications\rq\rq\,project 2022E2Z4AK – funded by European Union – Next Generation EU  within the PRIN 2022 program (D.D. 104 - 02/02/2022 Ministero dell’Università e della Ricerca).
\begin{figure}[H]
		\centering
		\includegraphics[scale=0.35]{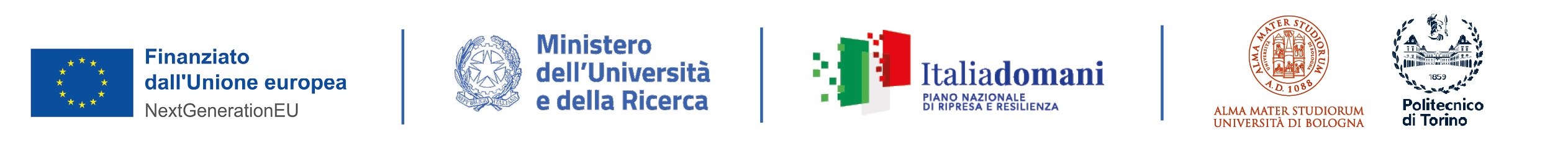}
			\end{figure}


\bigskip
\begin{thebibliography} {biblio} 

\bibitem[1] {BBCG} E.Ballico, A. Bernardi, M.Christandl , F.Gesmundo,  {\it On the partially symmetric rank of tensor products of W-states and other symmetric tensors. }   Atti della Accademia Nazionale dei Lincei, Classe di Scienze, Rendiconti Lincei Matematica E Applicazioni 30 (2018) DOI:10.4171/RLM/837, available at https://arxiv.org/abs/1803.01623 
\bibitem[2]{BF} E. Ballico, C.Fontanari,  {\it On the secant varieties to the osculating variety of a Veronese surface}, Central Europ. J. of Math. 1 (2003), 315-326.
\bibitem[3]{BFZ} W.Bruzda, S. Friedland, K. \.Zyczkowski, {\it Rank of a tensor and quantum entanglement.} Linear and Multilinear Algebra (2023), 1–64. https://doi.org/10.1080/03081087.2023.2211717
\bibitem[4]{BCCGO} A.Bernardi, E. Carlini, M.V.Catalisano, A. Gimigliano, A. Oneto  {\it The hitchhiker guide to secant varieties and tensor decomposition}. Mathematics 6 (2018,  special issue: "Decomposibility of Tensors", L.Chiantini Ed.);  DOI: 10.3390/math6120314. Also available at    http://www.mdpi.com/2227-7390/6/12/314/pdf
\bibitem[5] {BCGI} A.Bernardi, M.V.Catalisano, A.Gimigliano, M.Idà, {\it Osculating Varieties of Veronesean and their higher secant varieties.  } Canadian J. of  Math.  59, (2007), 488-502. 
\bibitem[6] {BCGI2} A.Bernardi, M.V.Catalisano, A.Gimigliano, M.Idà, {\it Secant varieties to Osculating Varieties of Veronese embeddings of  $\P^n$}.  J. of Algebra, 321, (2009),  982-1004
 (available also at:  http://arxiv.org/0807.2455v2)
\bibitem[7] {BGI} A.Bernardi, A.Gimigliano, M.Idà, {\it  Computing symmetric rank for symmetric tensors}. J. of Symbolic Computation, 46, (2011), 34–53 (available also at:  http://arxiv.org/abs/0908.1651)

\bibitem[BGM] {BGM} J. Briancon, M. Granger and P. Maisonobe, {\it Le nombre de modules du germe de courbe plane $x^a + y^b = 0$}. Math. Ann. 279, 535–551 (1988).
\bibitem[8] {CGI} S.Canino, A.Gimigliano, M.Idà, {\it On the Jacobian scheme of a plane curve.} Communications in Algebra, 1–11, 2024. https://doi.org/10.1080/00927872.2024.2384056. 
\bibitem[9] {CCG} E Carlini, M.V.Catalisano, A.V.Geramita, {\it The solution to the Waring problem for monomials and the sum of coprime monomials} J. of Symbolic Computation, 54, (2013), 9–35. 
\bibitem[11] {CGG} M.V.Catalisano, A.V.Geramita, A.Gimigliano, {\it Tensor rank, secant varieties to Segre varieties and fat points in multiprojective spaces}.  Lin. Alg. and Appl. 355, (2002) , 263-285. (see also the errata of the publisher: 367 (2003), 347-348). 
\bibitem[12] {CGG1} M.V.Catalisano, A.V.Geramita, A.Gimigliano, {\it On the secant varieties to the tangential varieties of a Veronesean}.   Proc. Am. Math. Soc.  130,  (2001), 975-985.
\bibitem[13] {CGG2} M.V.Catalisano, A.V.Geramita, A.Gimigliano, {\it  On the ideals of Secant Varieties to certain rational varieties. } J. of Algebra  319, (2008), 1913–1931  (also available at:  arXiv:math/0609054).
\bibitem[14] {CF} L. Chen and S. Friedland. {\it The tensor rank of tensor product of two three-qubit W-states is eight.} Linear Algebra and its Applications, 543, (2018), 1 – 16.
\bibitem[15] {COCOA} J.Abbott, A.M.Bigatti, L.Robbiano, {\it CoCoA: a system for doing Computations in Commutative Algebra.} Available at http://cocoa.dima.unige.it 
\bibitem[17]{FH} W.Fulton and J.Harris {\it Representation Theory, a first course}. Springer-Verlag, Berlin, Heidelberg, New York, 1999.
\bibitem[18]{F2} W.Fulton, {\it Intersection Theory}. Springer-Verlag, Berlin, Heidelberg, New York, 1998

\bibitem[19]{Ga} M. Ga\l azka. {\it Multigraded apolarity.}  Math. Nach., 296, (2023), 286–313. 
\bibitem[20]{Ge} A. Geramita, {\it Inverse systems of fat points: Waring’s problem, secant varieties of Veronese varieties, and parameter spaces for Gorenstein ideals}, in Queen’s Papers in Pure and Appl. Math. 102 (1996) 1-114.
\bibitem[G]{G}{J. M. Granger}, {\it Sur une espace de modules de germe de courbe plane.} Bull. Sci. Math. 26 s6r. 103,1979
\bibitem[21]{Ha}  B.Hartshorne, {\it Algebraic Geometry}. Springer, Grad. Texts in Math. 52, Berlin, Heidelberg, New York 1977.
\bibitem[22]{I} A. Iarrobino, {\it Inverse system of a symbolic power. III. Thin algebras and fat points}.Composito Math. 108 (1997) 319-356.
\bibitem[23]{J} J. Jelisiejew {\it Hilbert schemes of points and their applications} PhD dissertation, University of Warsaw (2017). https://arxiv.org/pdf/2205.10584.pdf
\bibitem[24]{LT} J.M.Landsberg, Z.Teitler, {\it On the ranks and border ranks of symmetric tensors}. Found. Comput. Math. 10 (2010), 339–366.
\bibitem[LP]{LP} O. A. Laudal, G. Pfister, {\it Local moduli and singularities}. Springer, Lecture notes in Mathematics 1310, New York/Berlin 1988.
\bibitem[26]{Pu} M.Pucci, {\it The Veronese Variety and Catalecticant Matrices}, J.of Algebra 202  (1998) 72-95.

\end{thebibliography}
\end{document}